\documentclass[titlepage,11pt]{article}
\usepackage{latexsym}
\usepackage{amsfonts}
\usepackage{amsmath}
\usepackage{tikz}
\usepackage{float}
\usepackage{lmodern}
\usepackage{underscore}
\usetikzlibrary{decorations.pathmorphing}

\newcommand\blackslug{\hbox{\hskip 1pt \vrule width 4pt height 8pt depth 1.5pt
        \hskip 1pt}}
\newcommand\bbox{\hfill \quad \blackslug \bigbreak}
\def\d{\hbox{-}}
\def\c{\hbox{-}\cdots\hbox{-}}
\def\ll{,\ldots,}

\title{Induced subgraphs of graphs with large chromatic number.
\\X. Holes of specific residue}
\author{Alex Scott\thanks{Supported by a Leverhulme Trust Research
Fellowship.}\\
Mathematical Institute, University of Oxford, Oxford OX2 6GG, UK
\\
\\
Paul Seymour\thanks{Supported by ONR grant N00014-14-1-0084 and NSF
grant DMS-1265563.}\\
Princeton University, Princeton, NJ 08544, USA}

\date{January 8, 2017; revised \today}

\newtheorem{thm}{}[section]

\newcommand{\Proof}{\noindent{\bf Proof.}\ \ }

\begin{document}
\maketitle
\begin{abstract}
A large body of research in graph theory concerns the induced subgraphs of graphs with large chromatic number, and especially which induced cycles must occur.
In this paper, we unify and substantially extend results from a number of previous papers, showing
that, for every positive integer $k$, every graph with  large chromatic number
contains either a large complete subgraph or induced cycles of all lengths 
modulo $k$.  
As an application, we prove two conjectures of Kalai and Meshulam from the 1990's connecting the chromatic number of a graph with the homology of its independence complex.
\end{abstract}

\section{Introduction}
All graphs in this paper are finite and have no loops or parallel edges. We denote the chromatic number of a graph $G$
by $\chi(G)$, and its clique number (the cardinality of its largest clique) by $\omega(G)$. A {\em hole} in $G$ means an induced subgraph
which is a cycle of length at least four. 

What can we say about the hole lengths in a graph $G$ with large chromatic number?  If $G$ is a complete  
graph then it has no holes at all, and the question becomes trivial.  
But if we bound the clique number of $G$ then the question becomes much more interesting, and much deeper.  
In an influential paper written thirty years ago,  Gy\'arf\'as~\cite{gyarfas} made a number of conjectures about 
induced subgraphs of graphs with large chromatic number and bounded clique number.  Three of these conjectures, concerning holes, are
particularly well-known:
\begin{thm}\label{gyarfasconj} For all $\kappa\ge 0$
\begin{itemize}
\item there exists $c$ such that every graph with chromatic number greater than $c$
contains either a complete subgraph on $\kappa$ vertices or a hole of odd length;
\item for all $\ell\ge 0$  there exists $c$ such that
every graph with chromatic number greater than $c$
contains either a complete subgraph on $\kappa$ vertices or a hole of length at least $\ell$;
\item for all $\ell\ge 0$ there exists $c$ such that
every graph with chromatic number greater than $c$
contains either a complete subgraph on $\kappa$ vertices or a hole 
whose length is odd and at least $\ell$.
\end{itemize}
\end{thm}

All three conjectures are now known to be true: the first was proved by the authors in~\cite{oddholes} (see \cite{cycles} for earlier work);
the second jointly with Maria Chudnovsky in~\cite{longholes}; and the third (which is a strengthening of the first two) jointly with Chudnovsky and Sophie Spirkl in~\cite{longoddholes}. 
 The analogous result for long even holes  is also known (it is enough to find two vertices joined by three long paths with no edges between them, and
this follows from results of \cite{bananas}).

Another intriguing result on holes was shown by Bonamy, Charbit and Thomass\'e~\cite{bonamy}, who
proved  a conjecture of Kalai and Meshulam by showing the following.

\begin{thm}\label{bonamy}
Every graph with sufficiently large chromatic number contains either a triangle or a hole of length $0$ modulo $3$.
\end{thm}

In this paper we prove the following theorem, which contains all the results mentioned above as special cases.

\begin{thm}\label{mainthm} 
For all $\kappa,\ell\ge 0$ there exists $c$ such that
every graph $G$ with $\chi(G)>c$ and $\omega(G)\le \kappa$
contains holes of every length modulo $\ell$.
\end{thm}

Note that this result allows us to demand a {\em long} hole of length $i$ modulo $j$ by 
taking $\ell=Nj$ for large $N$ and then choosing a suitable residue.  
Thus it implies all three Gy\'arf\'as conjectures; and it extends \ref{bonamy}
in several ways, 
allowing us to ask for any size of clique, and a hole of any residue and as long as we want.  
(Though we cannot demand a hole of any {\em specific} length: it is well-known that there are graphs
with arbitrarily large girth and chromatic number.)

We will in fact prove an even stronger statement.
We say $A,B\subseteq V(G)$ are {\em anticomplete} if $A\cap B=\emptyset$ and no vertex in $A$ has a neighbour in $B$;
and subgraphs $P,Q$ of $G$ are {\em anticomplete} if $V(P),V(Q)$ are anticomplete.  We prove the following.
\begin{thm}\label{superkalai}
Let $\kappa,n\ge 0$ be integers, and for $1\le i\le n$ let $p_i\ge 0$ and $q_i\ge 1$ be integers. Then there exists $c\ge 0$
with the following property. Let $G$ be a graph such that $\chi(G)>c$ and $\omega(G)\le \kappa$. Then
there are $n$ holes $H_1\ll H_{n}$ in $G$, pairwise anticomplete,
such that $H_i$ has length $p_i$ modulo $q_i$ for $1\le i\le n$.
\end{thm}

Let us restate this in slightly different language.
An ideal of graphs is {\em $\chi$-bounded} if there is a function $f$ such that $\chi(G)\le f(\omega(G))$ for every graph $G$ in the 
class. Thus \ref{superkalai} can be reformulated as:
\begin{thm}\label{boundedkalai}
Let $n\ge 0$ be an integer, and for $1\le i\le n$ let $p_i\ge 0$ and $q_i\ge 1$ be integers. Let $\mathcal{C}$
be the ideal of all graphs that do not contain $n$ pairwise anticomplete holes $H_1\ll H_n$ where
$H_i$ has length $p_i$ modulo $q_i$ for $1\le i\le n$.
Then $\mathcal{C}$ is $\chi$-bounded.
\end{thm}

\ref{boundedkalai} (or equivalently \ref{superkalai}) implies \ref{mainthm}, and also implies the main theorem of~\cite{complement}, which is the case of \ref{superkalai} when
$p_i=1$ and $q_i=2$ for each $i$. But it also has applications to further conjectures
of Kalai and Meshulam~\cite{kalai}, connecting graph theory with topology, and in particular with the homology of the independence complex of $G$. We discuss these in the final section.

Let us say a hole $H$ in $G$ is {\em $d$-peripheral} if $\chi(G[X])>d$, where $X$ is the set of vertices of $G$
that are not in $V(H)$ and have no neighbours in $V(H)$.
\ref{superkalai} follows easily from the following version of \ref{mainthm}, which will therefore be our main objective:
\begin{thm}\label{peripheral}
For all $\kappa,\ell,d\ge 0$ there exists $c$ such that
every graph $G$ with $\chi(G)>c$ and $\omega(G)\le \kappa$
contains 
$d$-peripheral holes of every length modulo $\ell$.
\end{thm}
\noindent{\bf Proof of \ref{superkalai}, assuming \ref{peripheral}.\ \ }
Let $\kappa, n$ and $p_i,q_i\;(1\le i\le n)$ be as in \ref{superkalai}. We may assume that
$n\ge 1$ and $\kappa\ge 2$, and we proceed by induction on $n$, for fixed $\kappa$.
Choose $d$
such that for every graph $G$ with $\chi(G)>d$ and $\omega(G)\le \kappa$,
there are $n-1$ holes $H_1\ll H_{n-1}$ in $G$, pairwise anticomplete,
where $H_i$ has length $p_i$ modulo $q_i$ for $1\le i\le n-1$. 
Let $c$ satisfy \ref{peripheral}
with $\ell$ replaced by $q_n$. We claim that $c$ satisfies \ref{superkalai}; for let 
$G$ be a graph such that $\chi(G)>c$ and $\omega(G)\le \kappa$. By \ref{superkalai}, $G$ has a $d$-peripheral 
hole $H_n$ of length $p_n$ modulo $q_n$. Let $X$ be the set of vertices of $G$ not in $H_n$ and with no neighbour in $H_n$.
Thus $\chi(G)>d$. From the inductive hypothesis, $G[X]$ has 
$n-1$ holes $H_1\ll H_{n-1}$ in $G$, pairwise anticomplete,
where $H_i$ has length $p_i$ modulo $q_i$ for $1\le i\le n-1$. But then $H_1\ll H_n$ satisfy the theorem.~\bbox

In this paper, we are also interested in holes of nearly
equal length.
In the triangle-free case, a result is known that is even stronger than \ref{mainthm}: we proved in~\cite{holeseq} that

\begin{thm}\label{holeseq} 
For all $\ell\ge 0$ there exists $c$ such that
every triangle-free graph with chromatic number greater than $c$
contains holes of $\ell$ consecutive lengths.
\end{thm}

We conjectured in~\cite{holeseq} that the same should be true if we exclude larger cliques:
\begin{thm}\label{moreholeseq}
{\bf Conjecture:} 
For all integers $\kappa, \ell\ge 0$,
there exists $c\ge 0$ such that 
every graph with chromatic number greater than $c$
contains either a complete subgraph on $\kappa$ vertices or 
holes of $\ell$ consecutive lengths.
\end{thm}
This conjecture remains open.  However, we make a small step towards it: we will show that under the same hypotheses, 
there are (long) holes of two consecutive lengths.
\begin{thm}\label{2holes}
For each $\kappa,\ell\ge 0$ there exists $c\ge 0$ such that 
every graph with chromatic number greater than $c$
contains either a complete subgraph on $\kappa$ vertices or 
holes of two consecutive lengths, both of length more than $\ell$.
\end{thm}
We have convinced ourselves that with a great deal of work, which we omit, we could get three consecutive ``long'' holes,
but so far that is the best we can do.

As in several other papers of this series, the proof of \ref{peripheral} examines whether there is an induced subgraph
of large chromatic number such that every ball of small radius in it has bounded chromatic number.
Let us make this more precise.
If $X\subseteq V(G)$, the subgraph of $G$ induced on $X$ is denoted by $G[X]$,
and we often write $\chi(X)$ for $\chi(G[X])$. The {\em distance} or {\em $G$-distance} between two vertices $u,v$
of $G$ is the length of a shortest path between $u,v$, or $\infty$ if there is no such path.
If $v\in V(G)$ and $\rho\ge 0$ is an integer, $N_G^{\rho}(v)$ or $N^{\rho}(v)$ denotes 
the set of all vertices $u$ with $G$-distance exactly
$\rho$  from $v$, and $N_G^{\rho}[v]$ or $N^{\rho}[v]$ denotes the set of all $u$ with $G$-distance at most $\rho$ from $v$.
If $G$ is a nonnull graph  and $\rho\ge 1$,
we define $\chi^{\rho}(G)$ to be the maximum of $\chi(N^{\rho}[v])$ taken over all vertices $v$ of $G$.
(For the null graph $G$ we define $\chi^{\rho}(G)=0$.)
Let $\mathbb{N}$ denote the set of nonnegative integers, and let $\phi:\mathbb{N}\rightarrow \mathbb{N}$ be a non-decreasing function.
For $\rho\ge 1$, let us say a graph $G$ is {\em $(\rho,\phi)$-controlled} if
$\chi(H)\le \phi(\chi^{\rho}(H))$ for every induced subgraph $H$ of $G$. Roughly, this says that in every induced subgraph $H$ of $G$ with
large chromatic number, there is a vertex $v$ such that $H[N^{\rho}_H[v]]$ has large chromatic number.
Let $\mathcal{C}$ be a class of graphs. We say $\mathcal{C}$ is an {\em ideal} if every induced subgraph of each member
of $\mathcal{C}$ also belongs to $\mathcal{C}$. If $\rho\ge 2$ is an integer, an ideal $\mathcal{C}$
is {\em $\rho$-controlled} if there is a nondecreasing function $\phi:\mathbb{N}\rightarrow \mathbb{N}$
such that every graph in $\mathcal{C}$ is $(\rho,\phi)$-controlled. For $\ell\ge 4$, an {\em $\ell$-hole} means a hole
of length exactly $\ell$. 
The proof of \ref{peripheral} breaks into two parts, the 2-controlled case and the $\rho$-controlled case when $\rho>2$
(because if we can be sure that all 2-balls have small chromatic number then it is easier to piece together paths to make holes
of any desired length.) We will prove the following two complementary results, 
which together imply \ref{peripheral}:

\begin{thm}\label{radthm}
Let $\rho\ge 2$ be an integer, and let $\mathcal{C}$ be a $\rho$-controlled ideal of graphs.
Let $\ell\ge 24$ if $\rho=2$, and $\ell\ge 8\rho^2+6\rho$ if $\rho>2$. Then for all $\kappa,d\ge 0$, 
there exists $c\ge 0$ such that every graph $G\in \mathcal{C}$ with $\omega(G)\le \kappa$ and $\chi(G)>c$ 
has a $d$-peripheral $\ell$-hole.
\end{thm}

\begin{thm}\label{uncontrolled}
For all integers $\ell\ge 2$ and $\tau,d\ge 0$ there is an integer $c\ge 0$ with the following property.
Let $G$ be a graph such that $\chi^8(G)\le \tau$, and every induced subgraph $J$ of $G$ with
$\omega(J)<\omega(G)$
has chromatic number at most $\tau$.
If $\chi(G)>c$ then there are $\ell$ $d$-peripheral holes in $G$ with lengths of all possible values modulo $\ell$.
\end{thm}
\noindent {\bf Proof of \ref{peripheral}, assuming \ref{radthm} and \ref{uncontrolled}.\ \ }
Let $\kappa,\ell,d\ge 0$, and let $\mathcal{C}$ be the ideal of graphs with clique number at most $\kappa$ and with no $d$-peripheral
hole of some 
length modulo $\ell$. 
By \ref{uncontrolled}, for each $\tau\ge 0$ there exists $c_{\tau}$ such that every $G\in \mathcal{C}$ with $\chi^8(G)\le \tau$ 
satisfies $\chi(G)\le c_{\tau}$,
and so $\mathcal{C}$ is 8-controlled. By \ref{radthm} the theorem follows. This proves \ref{peripheral}.~\bbox

We prove
the 2-controlled case of \ref{radthm} in the next section, and the $\rho>2$ case in section 3, deducing \ref{radthm} at the end of section 3. 
We prove \ref{uncontrolled} in section~\ref{sec:uncontrolled}, completing the proof of \ref{peripheral}; 
and 
prove the theorem about two consecutive long holes in section \ref{sec:2holes}.

\section{2-control}
First we handle the 2-controlled case. The proof here is very much like part of the proof of 
theorem 4.8 of~\cite{longoddholes}; the main difference is a strengthening of theorem 4.5 of that paper.
First we need some definitions.
If $G$ is a graph and $B,C\subseteq V(G)$, we say that $B$ {\em covers} $C$ if $B\cap C=\emptyset$ and every vertex in $C$
has a neighbour in $B$.
Let $G$ be a graph, let $x\in V(G)$, let $N$ be some set of neighbours of $x$, and let $C\subseteq V(G)$
be disjoint from $N\cup \{x\}$, such that $x$ is anticomplete to $C$ and $N$ covers $C$. In this situation
we call $(x,N)$ a {\em cover} of $C$ in $G$. For $C,X\subseteq V(G)$, a {\em multicover of $C$} in $G$
is a family $(N_x:x\in X)$
such that
\begin{itemize}
\item $X$ is stable;
\item for each $x\in X$, the pair $(x,N_x)$ is a cover of $C$;
\item for all distinct $x,x'\in X$, the vertex $x'$ is anticomplete to $N_x$ (and in particular all the sets
$\{x\}\cup N_x$ are pairwise disjoint).
\end{itemize}
Its {\em length} is $|X|$, and the multicover $(N_x:x\in X)$ is {\em stable} if each of the sets $N_x\;(x\in X)$ is stable.
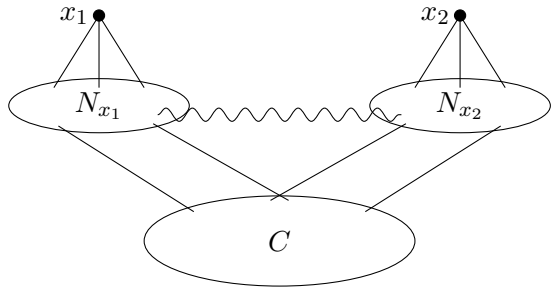
\begin{figure}[H]
\centering
\tikzset{snake it/.style={decorate, decoration=snake}}
\begin{tikzpicture}[scale=.6,auto=left]
\node at (0,-5) {$C$};
\node [left] at (-4,0) {$x_1$};
\node [left] at (4,0) {$x_2$};
\node at (-4,-2) {$N_{x_1}$};
\node at (4,-2) {$N_{x_2}$};
\tikzstyle{every node}=[inner sep=1.5pt, fill=black,circle,draw]
\node (x1) at (-4,0) {};
\node (x2) at (4,0) {};
\draw (-4,-2) ellipse (2cm and .6cm);
\draw (4,-2) ellipse (2cm and .6cm);
\draw (0,-5) ellipse (3cm and 1cm);
\draw (-4,0) -- (-4,-1.6);
\draw (-4,0) -- (-3,-1.6);
\draw (-4,0) -- (-5,-1.6);
\draw (4,0) -- (4,-1.6);
\draw (4,0) -- (3,-1.6);
\draw (4,0) -- (5,-1.6);
\draw (-2.8, -2.4) -- (0.2,-4.1);
\draw (-4.9, -2.45) -- (-1.9,-4.35);
\draw (2.8, -2.4) -- (-0.2,-4.1);
\draw (4.9, -2.45) -- (1.9,-4.35);
\path[draw=black, snake it] (-2.7, -2.2) -- (2.7, -2.2);

\end{tikzpicture}

\caption{A multicover of length two (the wiggle indicates possible edges)} \label{fig:1}
\end{figure}

Let $(N_x:x\in X)$ be a multicover of $C$, let $X'\subseteq X$, and for each $x\in X'$ let $N_x'\subseteq N_x$;
and let $C'\subseteq C$ be covered by each of the
sets $N_x'\;(x\in X')$. Then $(N_x':x\in X')$ is a multicover of $C'$, and we say it is {\em contained} in $(N_x:x\in X)$.

Again, let $(N_x:x\in X)$ be a multicover of $C$. Let $P$ be an induced path of $G$ with the following properties:
\begin{itemize}
\item $P$ has length three or five;
\item the ends of $P$ are in $X$;
\item no vertex of $X$ not an end of $P$ belongs to or has a neighbour in $V(P)$; and
\item every vertex of $P$ belongs to $X\cup \bigcup_{x\in X}N_x\cup C$.
\end{itemize}
Let us call such a path $P$ an {\em oddity} for the multicover. 

\begin{figure}[H]
\centering
\tikzset{snake it/.style={decorate, decoration=snake}}
\begin{tikzpicture}[scale=.6,auto=left]
\tikzstyle{every node}=[inner sep=1.5pt, fill=black,circle,draw]
\node (x1) at (-12,0) {};
\node (x2) at (-4,0) {};

\draw (-12,-2) ellipse (2cm and .6cm);
\draw (-4,-2) ellipse (2cm and .6cm);
\draw (-8,-5) ellipse (3cm and 1cm);
\draw (-12,0) -- (-12,-1.6);
\draw (-12,0) -- (-11,-1.6);
\draw (-12,0) -- (-13,-1.6);
\draw (-4,0) -- (-4,-1.6);
\draw (-4,0) -- (-5,-1.6);
\draw (-4,0) -- (-3,-1.6);
\draw (-10.8, -2.4) -- (-7.8,-4.1);
\draw (-12.9, -2.45) -- (-9.9,-4.35);
\draw (-5.2, -2.4) -- (-8.2,-4.1);
\draw (-3.1, -2.45) -- (-6.1,-4.35);

\node (y1) at (4,0) {};
\node (y2) at (12,0) {};
\draw (4,-2) ellipse (2cm and .6cm);
\draw (12,-2) ellipse (2cm and .6cm);
\draw (8,-5) ellipse (3cm and 1cm);
\draw (4,0) -- (4,-1.6);
\draw (4,0) -- (5,-1.6);
\draw (4,0) -- (3,-1.6);
\draw (12,0) -- (12,-1.6);
\draw (12,0) -- (11,-1.6);
\draw (12,0) -- (13,-1.6);
\draw (5.2, -2.4) -- (8.2,-4.1);
\draw (3.1, -2.45) -- (6.1,-4.35);
\draw (10.8, -2.4) -- (7.8,-4.1);
\draw (12.9, -2.45) -- (9.9,-4.35);

\node (a1) at (-11.2,-2) {};
\node (a2) at (-4.8,-2) {};
\node (b1) at (4.5,-2) {};
\node (b2) at (11.5,-2) {};
\node (c1) at (7,-4.5) {};
\node (c2) at (9,-4.5) {};

\foreach \from/\to in {x1/a1,x2/a2,a1/a2,y1/b1,y2/b2,b1/c1,b2/c2,c1/c2}
\draw [very thick] (\from) -- (\to);

\end{tikzpicture}

\caption{Oddities} \label{fig:2}
\end{figure}
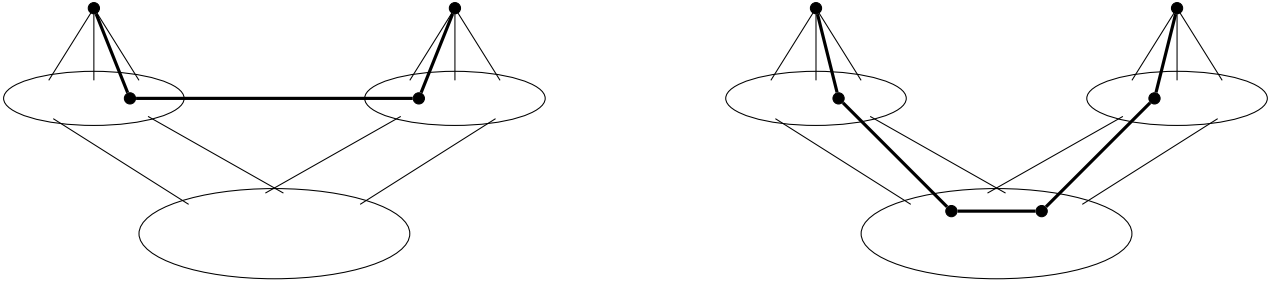

If $(N_x:x\in X)$ is a multicover of $C$, with an oddity $P$, and
$(N_x':x\in X')$ is a multicover of $C'\subseteq C$ contained in $(N_x:x\in X)$, and $V(P)$ is anticomplete
to $X'\cup \bigcup_{x\in X'}N_x'\cup C'$, we say that $(N_x':x\in X')$ is a multicover of $C'$ {\em compatible} with $P$.
Let $H$ be the subgraph induced on $\bigcup_{x\in X}N_x$; we call the clique number of $H$
the {\em cover clique number} of $(N_x:x\in X)$.

First we need to show the following:

\begin{thm}\label{getoddity}
Let $\tau,\kappa, m',c'\ge 0$ be integers, and let $0\le \kappa'\le \kappa$ be an integer.
Then there exist integers $m,c\ge 0$ with the following property.
Let $G$ be a graph such that
\begin{itemize}
\item $\omega(G)\le \kappa$;
\item $\chi(H)\le \tau$ for every induced subgraph $H$ of $G$
with $\omega(H)<\kappa$; and 
\item $G$ admits a
stable multicover $(N_x:x\in X)$ with length $m$, of a set $C$ with $\chi(C)>c$, with 
cover clique number at most $\kappa'$.
\end{itemize}
Then there is an oddity $P$ for the multicover, and a multicover $(N_x':x\in X')$
of $C'\subseteq C$ contained in $(N_x:x\in X)$ and compatible with $P$, such that $|X'|=m'$ and $\chi(C')>c'$.
\end{thm}
We proceed by induction on $\kappa'$, with $\tau,\kappa,m',c'$ fixed.
Thus, inductively, there exist $m_0,c_0\ge 0$ such that
the theorem holds if $m,c$ are replaced by $m_0,c_0$ respectively, and $\kappa'$ is replaced by any $\kappa_0$ with
$0\le \kappa_0<\kappa'$.
(Note that possibly $\kappa'=0$, when this statement is vacuous; in that case take $m_0=c_0=0$.)

Let $m=4+4m_0+2m'$. Define $c_m=4\tau+2^{m}(c_0+c')$, and
for $i=m-1\ll 1$ let $c_{i}=2c_{i+1}+\tau$. Let $c=2c_1+\tau$; we will show that $m,c$ satisfy the theorem.

Let $G$, $(N_x:x\in X)$ and $C$ be as in the theorem, where $|X|=m$, $\chi(C)>c$ and the cover clique
number of $(N_x:x\in X)$ is at most $\kappa'$. We may assume (because otherwise the theorem follows from the inductive hypothesis) that:
\\
\\
(1) {\em There is no multicover $(N_x':x\in X')$ of $C'\subseteq C$ contained in $(N_x:x\in X)$ with cover clique number less than $\kappa'$,
and with $|X'|=m_0$ and $\chi(C')>c_0$.}

\bigskip

Let $X=\{x_1\ll x_m\}$, and let us write $N_i$ for $N_{x_i}$ for $1\le i\le m$.
\\
\\
(2) {\em For $1\le i\le m$, there exist disjoint $C_i,D_i\subseteq C$ with $\chi(C_i), \chi(D_i)>c_i$, and $A_h\subseteq N_h$
for $1\le h\le i$, such that each $A_h$ covers one of $C_i, D_i$ and is anticomplete to the other.}
\\
\\
If $A\subseteq N_1\cup \cdots\cup N_m$, let $f(A)$ denote the set of vertices in $C$ with a neighbour in $A$.
Since $\chi(C)>c$, there exists $A_1\subseteq N_1$ minimal such that $f(A_1)$
has chromatic number more than $c_1$. Let $C_1=f(A_1)$ and $D_1=C\setminus C_1$.
From the minimality of $A_1$, it follows that $\chi(C_1)\le c_1+\tau$. Consequently
$\chi(D_1)>\chi(C)-(c_1+\tau)\ge c_1$. Thus (2) holds for $i = 1$. 
Now we assume that $i>1$ and $C_{i-1}, D_{i-1}$ and the sets $A_1\ll A_{i-1}$
satisfy (2) for $i-1$.
Choose $A_i\subseteq N_i$
minimal such that one of $\chi(f(A_i)\cap C_{i-1})$, $\chi(f(A_i)\cap D_{i-1})$ is more than $c_i$; say
the first (without loss of generality). Let $C_i=f(A_i)\cap C_{i-1}$. Now $\chi(f(A_i)\cap D_{i-1})\le c_i+\tau$, from the
minimality of $A_i$, so $\chi(D_i)>c_{i-1}-c_i-\tau\ge c_i$, where $D_i = D_{i-1}\setminus f(A_i)$. Thus $A_i$ covers $C_i$
and is anticomplete to $D_i$. This proves (2).

\bigskip

From (2) with $i = m$, each $A_i$ covers one of $C_m, D_m$ and is anticomplete to the other. By exchanging $C_m,D_m$
if necessary, we may assume that for at least $m/2$ values of $i$, $A_i$ covers $C_m$ and is anticomplete to $D_m$.
We may assume (by reordering $x_1\ll x_m$) that $A_i$ covers $C_m$ and is anticomplete to $D_m$ for $1\le i\le m/2$.
Let $B_i=N_i\setminus A_i$ for $1\le i\le m/2$.
\\
\\
(3) {\em There is an oddity $P$ for $(N_x:x\in X)$ with ends $x_1,x_2$ and with interior in
$B_1\cup B_2\cup D_m.$}
\\
\\
Since $\chi(D_m)>c_m\ge \tau$, there is a clique $Z\subseteq D_m$
with $|Z|=\kappa$. Now $N_1$ covers $D_m$, but $A_1$ is anticomplete to $D_m$, so $B_1$ covers $D_m$.
Similarly $B_2$ covers $D_m$. Choose a vertex $y_1\in B_{1}\cup B_{2}$
with as many neighbours in $Z$ as possible; and we may assume that
$y_1\in B_{1}$. Not every vertex of $Z$ is incident with $y_1$ since $\omega(G)\le \kappa$;
let $z_2\in Z$ be nonadjacent to
$y_1$. Choose $y_2\in B_{2}$ adjacent to $z_2$.
From the choice of $y_1$, there exists
$z_1\in Z$ adjacent to $y_1$ and not to $y_2$. If $y_1,y_2$ are nonadjacent, then $x_1\d y_1\d z_1\d z_2\d y_2\d x_2$ is an oddity,
and if $y_1,y_2$ are adjacent then $x_1\d y_1\d y_2\d x_2$ is an oddity. This proves (3).

\bigskip
Now there are at most four vertices of $P$ that have neighbours in $C_m$, and so there exists $F\subseteq C_m$ with
$\chi(F)>c_m-4\tau=2^{m}(c_0+c')$ that is anticomplete to $V(P)$. There are two vertices of $P$ in $N_1\cup N_2$, and those are the only 
vertices of $P$ that might have neighbours in $A_i$ for $3\le i\le m/2$. Let these vertices be $p,q$, and for $3\le i\le m/2$
let $P_i$ be the set of vertices in $A_i$ adjacent to $p$, and $Q_i$ the set adjacent to $q$.

For each $v\in F$, let $I(v)$ be the set of $i$ with $3\le i\le m/2$
such that $v$ has a neighbour in $P_i$. For each subset $I\subseteq \{3\ll m/2\}$ with $|I|=m_0$, the chromatic number of
the set of $v\in F$ with $I\subseteq I(v)$ is at most $c_0$, by (1). Since there are at most $2^{m-1}$ such subsets $I$,
the set of vertices $v\in F$ with $|I(f)|\ge m_0$ has chromatic number at most $2^{m-1}c_0$; and similarly the set of vertices
adjacent to neighbours of $q$ in at least $m_0$ sets $A_i$ has chromatic number at most $2^{m-1}c_0$. Consequently
there exists $F'\subseteq F$ with
$$\chi(F')\ge \chi(F)-2^{m}c_0> 2^{m}c'$$
such that for each $v\in F'$, there are at most $2m_0$ values of $i\in \{3\ll m/2\}$ such that $v$ is adjacent to a neighbour
of $p$ or $q$ in $A_i$. There are only at most $2^{m}$ possibities for the set of these values, so there exists
$C'\subseteq F'$ with $\chi(C')\ge \chi(F')2^{-m}>c'$ such that all vertices in $C'$ have the same set of values,
and in particular there exists $I\subseteq \{3\ll m/2\}$ with $|I|=m/2-2-2m_0=m'$ such that no vertex in $C'$
has a neighbour adjacent to $p$ or $q$ in any $A_i(i\in I)$. For each $i\in I$, let $N_{x_i}'$ be the set of vertices in $A_i$
nonadjacent to both $p,q$. Then $(N_{x}':x\in \{x_i:i\in I\})$ is a multicover of $C'$, contained in $(N_x:x\in X)$,
and compatible with $P$. This proves \ref{getoddity}.~\bbox

By three successive applications of \ref{getoddity} (one for each oddity), we deduce:

\begin{thm}\label{gettriple}
For all integers $\tau,\kappa\ge 0$, there exist integers $m,c\ge 0$ with the following property.
Let $G$ be a graph such that
\begin{itemize}
\item $\omega(G)\le \kappa$;
\item $\chi(H)\le \tau$ for every induced subgraph $H$ of $G$
with $\omega(H)<\kappa$; and
\item $G$ admits a
stable multicover $(N_x:x\in X)$ of a set $C$,
where $|X|=m$ and $\chi(C)>c$.
\end{itemize}
Then there are three oddities $P_1,P_2,P_3$ for the multicover, where $V(P_1), V(P_2), V(P_3)$
are pairwise anticomplete.
\end{thm}
(The same is true with ``three'' replaced by any other positive integer, but we only need three.)
Next we need:

\begin{thm}\label{findhole}
Let $\ell\ge 24$ be an integer.
Take the complete bipartite graph $K_{\ell,\ell}$, with bipartition $A,B$. Add three more edges joining three disjoint
pairs of vertices in $A$. Now subdivide every edge between $A$ and $B$ once, and subdivide each of the three additional
edges either
two or four times. The graph we produce has a hole of length $\ell$.
\end{thm}
We leave the proof to the reader (use the fact that if
$x,y,z\in \{3,5\}$ then $\ell$ is expressible as a sum of some or none of $x,y,z$ and at least three 4's).

A multicover $(N_x:x\in X)$ of $C$ is said to be {\em stably $k$-crested} if there are vertices $a_1\ll a_k$ and
vertices $a_{ix}\;(1\le i\le k, x\in X)$ of $G$,
all distinct, with the following properties:
\begin{itemize}
\item $a_1\ll a_k$ and the vertices $a_{ix}\;(1\le i\le k, x\in X)$ do not belong to $X\cup C\cup \bigcup_{x\in X}N_x$;
\item for $1\le i\le k$  and each $x\in X$, $a_{ix}$ is adjacent to $x$, and there are no other edges between the sets
$\{a_1\ll a_k\}\cup \{a_{ix}:1\le i\le k, x\in X\}$ and $X\cup C\cup \bigcup_{x\in X}N_x$;
\item for $1\le i\le k$  and each $x\in X$, $a_{ix}$ is adjacent to $a_i$, and there are no other edges between $\{a_1\ll a_k\}$ and $\{a_{ix}:1\le i\le k, x\in
X\}$
\item $a_1\ll a_k$ are pairwise nonadjacent;
\item for all $i,j\in \{1\ll k\}$ and all distinct $x,y\in X$, $a_{ix}$ is nonadjacent to $a_{jy}$.
\end{itemize}
(Thus the ``crest'' part is obtained from
$K_{k,|X|}$ by subdividing every edge once.)
We deduce:

\begin{thm}\label{crested}
Let $\ell\ge 24$, and let $\tau,\kappa\ge 0$. Then there exist $m,c\ge 0$ with the following property.
Let $G$ be a graph such that
\begin{itemize}
\item $\omega(G)\le \kappa$;
\item $\chi(J)\le \tau$ for every induced subgraph $H$ of $G$
with $\omega(H)<\kappa$; and
\item $G$ admits a stably
$\ell$-crested stable multicover $(N_x:x\in X)$ of a set $C$,
where $|X|=m$ and $\chi(C)>c$.
\end{itemize}
Then $G$ has a hole of length $\ell$.
\end{thm}
\Proof
Let $m,c$ satisfy \ref{gettriple}, choosing $m\ge \ell$ (note that if $m,c$ satisfy \ref{gettriple} then so do $m',c$ 
for $m'\ge m$.) By \ref{gettriple}, there are three oddities, pairwise anticomplete;
and the result follows from \ref{findhole}. This proves \ref{crested}.~\bbox

Theorem 4.4 of~\cite{longoddholes} says:
\begin{thm}\label{getbigtick}
For all $m,c,k,\kappa,\tau\ge 0$ there exist $m',c'\ge 0$ with the following property.
Let $G$ be a graph with $\omega(G)\le \kappa$, such that $\chi(H)\le \tau$ for every induced subgraph $H$ of $G$ with $\omega(H)<\kappa$.
Let $(N_x':x\in X')$ be a  multicover in $G$ of some set $C'$, such that $|X'|\ge m'$ and $\chi(C')> c'$.
Then there exist $X\subseteq X'$ with $|X|\ge m$, and $C\subseteq C'$ with $\chi(C)> c$, and a stable multicover $(N_x:x\in X)$
of $C$ contained in $(N_x':x\in X')$ that is stably $k$-crested.
\end{thm}

Combining \ref{crested} and \ref{getbigtick}, we deduce the following (a strengthening of 
theorem 4.5 of~\cite{longoddholes}):

\begin{thm}\label{nomult}
Let $\ell\ge 24$, and let $\tau,\kappa\ge 0$. Then there exist $m,c\ge 0$ with the following property.
Let $G$ be a graph such that
\begin{itemize}
\item $\omega(G)\le \kappa$;
\item $\chi(H)\le \tau$ for every induced subgraph $H$ of $G$
with $\omega(H)<\kappa$; and
\item $G$ admits a 
multicover with length $m$, of a set $C$ with $\chi(C)>c$.
\end{itemize}
Then $G$ has a hole of length $\ell$.
\end{thm}

We need the following, a consequence of theorem 9.7 of~\cite{chandeliers}. That involves ``trees of lamps'', but we
do not need to define
those here; all we need is that a cycle of length $\ell$ is a tree of lamps. (Note that what we call a
``multicover'' here is called a ``strongly-independent 2-multicover'' in that paper, and indexed in a slightly different way.)
\begin{thm}\label{findchand5}
Let $m,\kappa,c',\ell\ge 0$, and let $\mathcal{C}$ be a 2-controlled ideal,
such that for every $G\in \mathcal{C}$:
\begin{itemize}
\item $\omega(G)\le \kappa$;
\item $G$ does not admit a multicover of length $m$ of a set with chromatic number more than $c'$; and
\item $G$ has no hole of length $\ell$.
\end{itemize}
Then there exists $c$ such that all graphs in $\mathcal{C}$ have chromatic number at most $c$.
\end{thm}

\bigskip

Now we prove the main result of this section, that is, \ref{radthm} with $\rho=2$. 

\begin{thm}\label{rad2}
Let $\ell\ge 24$ and let $\mathcal{C}$ be a 2-controlled ideal of graphs. For all $\kappa,d\ge 0$ there exists $c$
such that every graph in $\mathcal{C}$ with clique number at most $\kappa$ and chromatic number more than $c$
has a $d$-peripheral hole of length $\ell$.
\end{thm}
\Proof
We proceed by induction on $\kappa$. The result holds for $\kappa\le 1$, so we assume that $\kappa\ge 2$ and
every graph in $\mathcal{C}$ with clique number less than $\kappa$ has chromatic number at most $\tau$.
Let $\mathcal{C}'$ be the ideal of graphs $G\in \mathcal{C}$ such that $\omega(G)\le \kappa$ and $G$ has no hole of 
length $\ell$. Choose $m,c'$ to satisfy \ref{nomult} (with $c$ replaced by $c'$). Choose $c''$ to satisfy
\ref{findchand5} (with $\mathcal{C}$ replaced by $\mathcal{C}'$, and $c$ replaced by $c''$). Let $c=\max(c'',d+\ell \tau)$; we claim that
$c$ satisfies the theorem. For let $G\in \mathcal{C}$ with $\omega(G)\le \kappa$ and $\chi(G)>c$. 
Suppose first that $G\in \mathcal{C}'$. By \ref{nomult},        
$G$ does not admit a multicover with length $m$ of a set with chromatic number more than $c'$. From
\ref{findchand5}, $\chi(G)\le c''$, a contradiction.

Thus $G\notin \mathcal{C}'$, and so $G$
has an $\ell$-hole $H$. For each vertex of $H$, its set of neighbours has chromatic number at most $\tau$; and so
the set of all vertices of $G$ that belong to or have a neighbour in $H$ has chromatic number at most $\ell\tau$. Since
$\chi(G)>d+\ell \tau$, it follows that $H$ is $d$-peripheral. This proves \ref{rad2}.~\bbox

\section{The $\rho$-controlled case for $\rho\ge 3$.}

Let $G$ be a graph. We say a {\em grading} of $G$ is a sequence $(W_1\ll W_n)$ of subsets of $V(G)$, pairwise
disjoint and with union $V(G)$. If $w\ge 0$ is such that $\chi(G[W_i])\le \tau$ for $1\le i\le n$
we say the grading is {\em $\tau$-colourable}. We say that $u\in V(G)$ is {\em earlier} than $v\in V(G)$
(with respect to some grading $(W_1\ll W_n)$) if $u\in W_i$ and $v\in W_j$ where $i<j$.

Let $G$ be a graph, and let $B,C\subseteq V(G)$, where $B$ covers $C$. Let $B=\{b_1\ll b_m\}$.
For $1\le i<j\le m$ we say that $b_i$ is {\em earlier} than $b_j$
(with respect to the enumeration $(b_1\ll b_m)$). For $v\in C$, let $i\in \{1\ll m\}$ be minimum such that $b_i,v$ are adjacent;
we call $b_i$ the {\em earliest parent} of $v$. An edge $uv$ of $G[C]$ is said to be {\em square} (with respect to the
enumeration
$(b_1\ll b_m)$) if
the earliest parent of $u$ is nonadjacent to $v$, and the earliest parent of $v$ is nonadjacent to $u$.
Let $B=\{b_1\ll b_m\}$, and let $(W_1\ll W_n)$ be a grading of $G[C]$. We say the enumeration $(b_1\ll b_m)$  of $B$ and the
grading $(W_1\ll W_n)$ are {\em compatible} if for all $u,v\in C$ with $u$ earlier than $v$, the earliest parent of $u$ is
earlier than
the earliest parent of $v$.

A graph $H$ is a {\em $\rho$-ball} if either $V(H)=\emptyset$ or there is a vertex $z\in V(H)$ such that every vertex of $H$
has $H$-distance at most $\rho$ from $z$; and we call $z$ a {\em centre} of the $\rho$-ball.
If $G$ is a graph, a subset $X\subseteq V(G)$ is said to be a {\em $\rho$-ball} if $G[X]$ is a
$\rho$-ball. (Note that there may be vertices of $G$ not in $X$ that have $G$-distance at most $\rho$ from $z$; and also, 
for a pair of vertices in $X$, their $G$-distance and their $G[X]$-distance may be different.)

\begin{thm}\label{greentouchrad}
Let $\phi$ be a nondecreasing function and $\rho\ge 3$, and let $G$ be a $(\rho,\phi)$-controlled graph. Let $\tau\ge 0$
such that $\chi^{\rho-1}(G)\le \tau$ and $\chi(J)\le \tau$ for every induced subgraph $J$ of $G$ with $\omega(J)<\omega(G)$.
Let $c\ge 0$ and let $(W_1\ll W_n)$ be a $\tau$-colourable grading of $G$. Let $H$ be a subgraph of $G$ (not necessarily induced)
with $\chi(H)>\tau+1+\phi(c+\tau)$, and such that $W_i\cap V(H)$ is stable in $H$ for each $i\in \{1\ll n\}$. 
Then there is an edge $uv$ of $H$, and a $\rho$-ball $X$ of $G$,
such that
\begin{itemize}
\item $u,v$ are both earlier than every vertex in $X$;
\item $v$ has a $G$-neighbour in $X$, and $u$ does not; and
\item $\chi(G[X])>c$.
\end{itemize}
\end{thm}
\Proof
Let us say that $v\in V(G)$ is {\em internally active} if there is a $\rho$-ball $X\ni v$ with $\chi(X)>c+\tau$ such that no vertex of $X$ is 
earlier than $v$. (Note that $X\cap W_i$ may have more than one element, so there may be vertices in $X$ that are neither
earlier nor later than $v$.) Let $R_1$ be the set of internally active vertices. We claim first:
\\
\\
(1) {\em $\chi(G\setminus R_1)\le \phi(c+\tau)$.}
\\
\\
For suppose not. Then since $G$ is $(\rho,\phi)$-controlled, there is a $\rho$-ball $X\subseteq V(G)\setminus R_1$
with $\chi(G)>c+\tau$, which therefore contains an internally active vertex, a contradiction. This proves (1).

\bigskip

Let us say $v\in V(G)$ is {\em externally active} if there is a $\rho$-ball $X$ of $G$ with $\chi(X)>c+\tau$ such that every vertex
of $X$ is later than $v$, and $v$ has an $H$-neighbour in $X$. Let $R_2$ be the set of externally active vertices. We claim:
\\
\\
(2) {\em $R_1\setminus R_2$ is stable in $H$.}
\\
\\
For suppose that $uv$ is an edge of $H$ with both ends in $R_1\setminus R_2$. Since each $W_i\cap V(H)$ is stable in $H$, we may assume that
$u$ is earlier than $v$. Since $v$ is internally active, there is a $\rho$-ball $X$ containing $v$ with $\chi(X)>c+\tau$ such that
no vertex of $X$ is earlier than $v$; but then $u$ is externally active, a contradiction. This proves (2).
\\
\\
(3) {\em There is a subset $Y\subseteq V(H)$ such that $H[Y]$ is connected and has chromatic number more than $\tau$,
and a $\rho$-ball $X$ of $G$ with $\chi(G[X])>c+\tau$, such that every vertex of $Y$ is earlier than every vertex of $X$, 
and some vertex of $Y$ has
a $H$-neighbour in $X$.}
\\
\\
Since $H$ has chromatic number more than $\tau+1+\phi(c+\tau)$, it follows from (1) and (2) that $\chi(H[R_2])>\tau$.
Let $Y$ be the vertex set of a component of $H[R_2]$ with maximum chromatic number. Choose $v\in Y$ such that no vertex
of $Y$ is later than $v$. Since $v$ is externally active, this proves (3).

\bigskip
Let $X,Y$ be as in (3). If some vertex of $Y$ has no $G$-neighbour in $X$, then since $H[Y]$ is connected, there is an
edge $uv$ of $H[Y]$ such that $v$ has a $G$-neighbour in $X$ and $u$ does not, and the theorem holds. We assume then
that every vertex of $Y$ has a $G$-neighbour in $X$. For each $y\in Y$, let $N(y)$ denote its set of $G$-neighbours in $X$.
Let $z$ be a centre of $X$, and for $0\le i\le \rho$ let $L_i$
be the set of vertices in $X$ with $G[X]$-distance $i$ to $z$. Thus $L_0\cup\cdots\cup L_{\rho}=X$.
Let $Y_0$ be the set of all $y\in Y$ with $N(y)\subseteq L_{\rho-1}\cup L_{\rho}$.
\\
\\
(4) {\em $Y_0\ne\emptyset$.}
\\
\\
Since $\chi(H[Y])>\tau$, it follows that $\chi(G[Y])>\tau$, and so some vertex $y\in Y$ has $G$-distance at least $\rho$
from $z$. Consequently $N(y)\subseteq L_{\rho-1}\cup L_{\rho}$. This proves (4).

\bigskip

Choose $y\in Y_0$, if possible with the additional property that $N(y)\cap L_{\rho-1}=\emptyset$.
Let $U$ be the set 
of vertices in $L_{\rho}$ with a neighbour in $N(y)\cap L_{\rho-1}$. 
\\
\\
(5) {\em There is a vertex $y'$ of $Y$ with $N(y')\not\subseteq N(y)\cup U$.}
\\
\\
For there is a vertex $y'\in Y$ with $G$-distance at least $\rho$ from $y$, since $\chi(G[Y])>\tau$. Since $\rho>2$,
$N(y)\cap N(y')=\emptyset$. If $N(y')\subseteq U$, then $y'\in Y_0$ and $N(y')\cap L_{\rho-1}=\emptyset$;
but then $N(y)\cap L_{\rho-1}=\emptyset$ from the choice of $y$, and so $U=\emptyset$, a contradiction. Thus 
$N(y')\not\subseteq U$. This proves (5).

\bigskip

Now $X\setminus (N(y)\cup U)$ is a $\rho$-ball $X'$ say, and some vertex (namely $y'$) of $Y$ has a $G$-neighbour in it,
and another (namely $y$) has no $G$-neighbour in it. Since $H[Y]$ is connected, there is an edge $uv$ of $H[Y]$
such that $v$ has a $G$-neighbour in $X'$ and $u$ does not. But $\chi(X)>c+\tau$, and
every vertex in $N(y)\cup U$ has $G$-distance at most two from $y$ and so $\chi(N(y)\cup U)\le \tau$,
and consequently $\chi(X')\ge \chi(X)-\tau>c$.
This proves \ref{greentouchrad}.~\bbox

We also need the following, proved in~\cite{longoddholes}:

\begin{thm}\label{getgreenedge}
Let $G$ be a graph, and let $B,C\subseteq V(G)$, where $B$ covers $C$. Let every induced subgraph $J$ of $G$ with
$\omega(J)<\omega(G)$
have chromatic number at most $\tau$.
Let the enumeration $(b_1\ll b_m)$ of $B$ and the grading $(W_1\ll W_n)$ of $G[C]$ be compatible. Let $H$ be the subgraph of $G$
with vertex set $C$ and edge set the set of all square edges. Let $(W_1\ll W_n)$ be $\tau$-colourable; then
$\chi(G[C])\le \tau^2\chi(H)$.
\end{thm}

We deduce:
\begin{thm}\label{greenedgerad}
Let $\phi$ be a nondecreasing function and $\rho\ge 3$, and let $G$ be a $(\rho,\phi)$-controlled graph. Let $\tau\ge 0$
such that $\chi^{\rho-1}(G)\le \tau$ and $\chi(J)\le \tau$ for every induced subgraph $J$ of $G$ with $\omega(J)<\omega(G)$.
Let $B,C\subseteq V(G)$, where $B$ covers $C$. 
Let the enumeration $(b_1\ll b_m)$ of $B$ and the grading $(W_1\ll W_n)$ of $G[C]$ be compatible.
Let $(W_1\ll W_n)$ be $\tau$-colourable, and let $\chi(G[C])>\tau^2 (\tau+1+\phi(c+\tau))$.
Then there is a square edge $uv$, and a $\rho$-ball $X$ of $G$,
such that
\begin{itemize}
\item $u,v$ are both earlier than every vertex in $X$;
\item $v$ has a neighbour in $X$, and $u$ does not; and
\item $\chi(X)>c$.
\end{itemize}
\end{thm}
\Proof Let $H$ be as in \ref{getgreenedge}. By \ref{getgreenedge}, 
$\chi(G[C])\le \tau^2\chi(H)$.
Since $\chi(G[C])>\tau^2 (\tau+1+\phi(c+\tau))$ and $\chi^1(G)\le \tau$, it follows that $\chi(H)>\tau+1+\phi(c+\tau)$.
By \ref{greentouchrad} applied to $G[C]$ and $H$, we deduce
that there is an edge $uv$ of $H$, and a $\rho$-ball $X$ of $G$, satisfying the theorem. This proves \ref{greenedgerad}.~\bbox

A {\em $\rho$-comet} $(\mathcal{P}, X)$ in a graph $G$ consists of a set $\mathcal{P}$ of induced paths, each with the same pair of
ends $x,y$ say, and a $\rho$-ball $X$, such that
$y$ has a neighbour in $X$ and no other vertex of any member of $\mathcal{P}$ has a neighbour in $X$. We call $x$ the {\em tip}
of the $\rho$-comet, and $\chi(X)$ its {\em chromatic number}, and the set of lengths of members of $\mathcal{P}$ its {\em spectrum}.

\begin{thm}\label{twotails}
Let $\phi$ be a nondecreasing function, and let $\rho\ge 3$ and $\tau\ge 0$.
For all integers $c\ge 1$ there exists $c'\ge 0$ with the following property.
Let $G$ be a $(\rho,\phi)$-controlled graph such that 
$\chi^{\rho-1}(G)\le \tau$ and $\chi(J)\le \tau$ for every induced subgraph $J$ of $G$ with $\omega(J)<\omega(G)$.
Let $x\in V(G)$, and let $V(G)\setminus \{x\}$ be a $\rho$-ball, such that $x$ has a neighbour in $G\setminus x$.
Let $\chi(V(G)\setminus \{x\})>c'$.
Then there is a $\rho$-comet $(\{P,Q\}, C)$ in $G$ with tip $x$ and chromatic number more than $c$, where
$|E(Q)|=|E(P)|+1$, and $|E(P)|\le 2\rho+1$.
\end{thm}
\Proof Let $c' = 2\tau^2 (\tau+1+\phi(c+\tau))$, and let $G,x$ be as in the theorem.
Since $V(G)\setminus \{x\}$ is a $\rho$-ball, every vertex of $G$ has $G$-distance at most $2\rho+1$ from $x$;
for $0\le k\le 2\rho+1$ let $L_k$ be the set of vertices of $G$ with $G$-distance exactly $k$
from $x$. Since $\chi(V(G)\setminus \{x\})>c'$, there exists $k$
such that $\chi(L_k)>c'/2$. Since $\chi^2(G)\le \tau$ it follows that $k\ge 3$.
Let $(b_1\ll b_n)$ be an enumeration of $L_{k-1}$, and for $1\le i\le n$ let $W_i$ be the set of vertices in $L_k$
that are adjacent to $b_i$ but not to $b_1\ll b_{i-1}$. Then $(W_1\ll W_n)$ is a $\tau$-colourable grading of 
$G[L_k]$, compatible with $(b_1\ll b_n)$. 

Since $\chi(L_k)>\tau^2 (\tau+1+\phi(c+\tau))$,
by \ref{greenedgerad} 
there is a square edge $uv$ of $G[L_k]$, and a $\rho$-ball $C$ of $G[L_k]$,
such that
\begin{itemize}
\item $u,v$ are both earlier than every vertex in $C$;
\item $v$ has a neighbour in $C$, and $u$ does not; and
\item $\chi(C)>c$.
\end{itemize}
Let $u',v'$ be the earliest parents of $u,v$ respectively.
Let $P$ consist of the union of the path $v\d v'$ and a path of length $k-1$ between $v',x$ with interior in $L_1\ll L_{k-2}$;
and let $Q$ consist of the union of the path $v\d u\d u'$ and a path of length $k-1$ between $u',x$ with interior in $L_1\ll L_{k-2}$.
Then $|E(Q)|=|E(P)|+1$ and $|E(P)|\le 2\rho$. Moreover, no vertex in $C$ has a neighbour in $P\cup Q$
different from $v$, since all vertices in $C$ are later than $u,v$. This proves \ref{twotails}.~\bbox

By repeated application of \ref{twotails} we deduce:

\begin{thm}\label{manytails}
Let $\phi$ be a nondecreasing function, and let $\rho\ge 3$, $\ell\ge \rho(8\rho+6)$ and $\tau\ge 0$.
For all integers $c\ge 1$ there exists $c'\ge 0$ with the following property.
Let $G$ be a $(\rho,\phi)$-controlled graph such that
$\chi^{\rho-1}(G)\le \tau$ and $\chi(J)\le \tau$ for every induced subgraph $J$ of $G$ with $\omega(J)<\omega(G)$.
Let $x\in V(G)$, and let $V(G)\setminus \{x\}$ be a $\rho$-ball, such that $x$ has a neighbour in $G\setminus x$.
Then there is a $\rho$-comet $(\mathcal{P}, X)$ in $G$ with tip $x$ and chromatic number more than $c$, such that
its spectrum includes $\{\ell+i:\:0\le i\le 2\rho+3\}$.
\end{thm}
\Proof Let $c_{\ell+1}=c$, and for $i=\ell\ll 1$ let \ref{twotails} be satisfied setting $c = c_{i+1}$ and $c'=c_{i}$.
Let $c'=c_1$.
\\
\\
(1) {\em For all $k\ge 1$ there exists $p_k$ with $1\le p_k\le 2\rho$ and a $\rho$-comet in $G$ with tip $x$,
chromatic number more than $c_k$, 
and spectrum including $\{p_1+\cdots+p_k+i:\:1\le i\le k\}$.}
\\
\\
By hypothesis there is a $\rho$-comet in $G$ with chromatic number more than $c_1$, tip $x$ and spectrum $\{1\}$, so the statement
holds when $k=1$, setting $p_1=0$; and it follows for $k\ge 2$ by repeated application of \ref{twotails}.
This proves (1).

\bigskip
Now $p_1\ll p_{\ell}$ exist and sum to at least $\ell$, so there exists $k\le \ell$ maximum such that 
$$p_1+\cdots+p_{k}\le \ell.$$
Since
$p_1+\cdots+p_{4\rho+3}< 2(4\rho+3)\rho\le \ell$, it follows that $k\ge 4\rho+3$.
From the maximality of $k$, and since $p_{k+1}\le 2\rho$, it follows that 
$p_1+\cdots+p_{k}> \ell-2\rho$. Consequently the spectrum of the corresponding $\rho$-comet contains 
$\{\ell+i:\:0\le i\le 2\rho+3\}$. This proves \ref{manytails}.~\bbox

\begin{thm}\label{manyholes}
Let $\phi$ be a nondecreasing function, and let $\rho\ge 3$, $\ell\ge 8\rho^2+6\rho$, and $d,\tau\ge 0$.
Then there exists $c$ with the following property.
Let $G$ be a $(\rho,\phi)$-controlled graph with $\chi(G)>c$ such that
$\chi^{\rho-1}(G)\le \tau$ and $\chi(J)\le \tau$ for every induced subgraph $J$ of $G$ with $\omega(J)<\omega(G)$.
Then there is a $d$-peripheral $\ell$-hole in $G$.
\end{thm}
\Proof
Define $c_4=\ell(2\rho+4)\tau$. Choose $c_3$ such that \ref{manytails} is satisfied replacing $c,c', \ell$ 
by $c_4,c_3,\ell-6\rho$ respectively.
Let $c_2 =\rho\tau +\phi(c_3)$,  $c_1=\tau\phi(c_2)$, and let $c=\max(\phi(c_1),\ell\tau+d)$. Let $G$ be as in the theorem with $\chi(G)>c$.
Since $\chi(G)>c\ge \phi(c_1)$, there exists $z\in V(G)$ such that, denoting the set of vertices
of $G$ with $G$-distance $i$ from $z$ by $L_i$, we have $\chi(L_{\rho})>c_1$. Since $L_1$ is $\tau$-colourable, there is a stable
subset $A$ of $L_1$ such that the set $B$ of vertices in $L_{\rho}$ that are descendants of vertices in $A$ has chromatic number
more than $c_1/\tau=\phi(c_2)$. Consequently there is a $\rho$-ball $C\subseteq B$ with $\chi(C)>c_2$.
Choose $D\subseteq A$ minimal such that every vertex in $C$ has an ancestor in $D$. Let $v_1\in D$; then
there exists $v_{\rho-1}\in L_{\rho-1}$ with a neighbour in $C$ such that $v_1$ is its only ancestor in $D$. 
Let $v_1\d v_2\c v_{\rho-1}$ be an induced path, where
$v_i\in L_i$ for $1\le i\le {\rho}-1$. The set of vertices in $C$
with $G$-distance less than $\rho$ from one of $v_1,v_2\ll v_{\rho-1}$ has chromatic number at most $\rho\tau$, and so the set $E$ of
vertices in $C$ with $G$-distance at least $\rho$ from each of $v_1,v_2\ll v_{\rho-1}$ has chromatic number more than
$c_2 -\rho\tau =\phi(c_3)$. Consequently there is a $\rho$-ball $F\subseteq E$, with chromatic number more than $c_3$.

Since $C$ is a $\rho$-ball and $v_{\rho-1}$ has a neighbour in $C$, there is an induced path $P$ of $G[C\cup \{v_{\rho}\}]$
from $v_{\rho-1}$ to some vertex
$x\in C$ with a neighbour in $F$, of length at most $2\rho$, such that no vertex of $P$ different from $x$
has a neighbour in $F$. By \ref{manytails} applied to $x,F$, since $\chi(F)> c_3$, there is a vertex $v\in F$, $2\rho+4$ induced paths
$P_0\ll P_{2\rho+3}$ of $G[F\cup\{v_{\rho-1}\}]$ between $x,v$, and a $\rho$-ball $X\subseteq F$, such that:
\begin{itemize}
\item  $|E(P_i)|=\ell-6\rho+i$ for $0\le i\le 2\rho+3$;
\item $V(P_i)\cap X=\emptyset$ for $0\le i\le 2\rho+3$;
\item $v$ has a neighbour in $X$ and no other vertex of $P_i$
has a neighbour in $X$, for $0\le i\le 2\rho+2$; and
\item $\chi(X)>c_4$.
\end{itemize}
Now every vertex of $X$ has $G$-distance at least $\rho$ from each of $v_1\ll v_{\rho-1}$, but there may be vertices in $X$
with $G$-distance less than $\rho$ to a vertex in $P$ or in one of $P_0\ll P_{2\rho+3}$. The union of these paths has at most
$\ell(2\rho+4)$ vertices (in fact much fewer), and since $\chi(X)>\ell(2\rho+4)\tau$, there exists a vertex in $X$ with $G$-distance
at least $\rho$ from all vertices of these paths. Let $y\in L_{\rho-1}$ be adjacent to this vertex, and let $Q$ be an induced path
between $y,x$ with interior in $X$, of length at most $2\rho+2$. Let $R$ be a path of length $\rho-1$
between $y,z$ with interior in $D\cup L_2\cup L_3\cup \cdots\cup L_{\rho-2}$. 

The union of the four paths $z\d v_1\d v_2\c v_{\rho-1}, P,Q,R$ has length at most $6\rho$, and at least $4\rho-3$, since
$P,Q$ have lengths at least $\rho$ and at least $\rho-1$ respectively. Let their union have length $j$
where $4\rho-3\le j\le 6\rho$. Let $i = 6\rho-j$; then $0\le i\le 2\rho+3$, and so $P_i$ is defined and has length $\ell-6\rho+i+j$.
Consequently the union of the five paths
$z\d v_1\d v_2\c v_{\rho-1}$, $P$, $P_i$, $Q$ and $R$ is a cycle $H_i$ of length $\ell$, and we claim it is induced.
Certainly it is a cycle; suppose that it is 
not induced, and so there is an edge $ab$ say that joins two nonconsecutive vertices of $H_i$. It follows
that $a,b$ do not both belong to any of its five constituent paths. Certainly $a,b\ne z$.
Suppose that $a=v_i$ for some $i$. Since every vertex in $E$ has $G$-distance at least $\rho$ from $v_i$, it follows that $b\notin V(P_i)$
and $b\notin V(Q)$. Also $b\notin V(P)$ since every vertex of $V(P)\setminus \{v_{\rho-1}\}$ belongs to $L_k$ and $P$
is an induced path containing $v_{\rho-1}$. Thus $b\in R$. Since the $G$-distance between $y,a$ is at least $\rho-1$, and $R$
has length $\rho-1$, it follows that $b\in L_1$ and so $i \in \{1,2\}$; but $i\ne 1$ since $A$ is stable, and $i\ne 2$ since
$v_{\rho-1}$, and hence $v_2$, has a unique ancestor in $D$. This proves that $a,b\notin \{v_1\ll v_{\rho-1}\}$.

Next suppose that $a\in V(R)$, and so either $b\in V(P)\setminus \{v_{\rho-1}\}$, or $b\in V(P_i\cup Q\setminus \{y\})$. 
In either case $b\in L_k$, and so $a=y$; hence $b\notin V(Q)$ since $Q$ is induced and $y\in V(Q)$, and $a\notin V(P\cup P_i)$
since $y$ has a neighbour in $X$ with $G$-distance at least $\rho$ from each vertex of $V(P\cup P_i)$, and $\rho\ge 3$.
This proves that $a,b\notin V(R)$. Next suppose that $a\in V(P)\setminus \{x\}$. Then $b\in F$; but no vertex of $P$ 
except $x$ has a neighbour in $F$, a contradiction. Finally, suppose that $a\in V(P_i)$ and $b\in V(Q)$. No vertex of $P_i$
has a neighbour in $X$ except $v$, and $v\in V(Q)$, a contradiction. This proves that $H_i$ is an $\ell$-hole. Now the
set of vertices of $G$ that belong to or have a neighbour in $H_i$ has chromatic number at most $\ell\tau$,
and since $\chi(G)>c\ge \ell\tau+d$, it follows that $H_1$ is $d$-peripheral.
This proves \ref{manyholes}.~\bbox

Let us deduce \ref{radthm}, which we restate:
\begin{thm}\label{radthm2}
Let $\rho\ge 2$ be an integer, and let $\mathcal{C}$ be a $\rho$-controlled ideal of graphs.
Let $\ell\ge 24$ if $\rho=2$, and $\ell\ge 8\rho^2+6\rho$ if $\rho>2$. Then for all $\kappa,d\ge 0$,
there exists $c\ge 0$ such that every graph $G\in \mathcal{C}$ with $\omega(G)\le \kappa$ and $\chi(G)>c$
has a $d$-peripheral $\ell$-hole.
\end{thm}
\Proof 
By induction on $\kappa$ we may assume that there exists $\tau_1$ such that every graph in $\mathcal{C}$ with clique number
less than $\kappa$ and no $d$-peripheral $\ell$-hole has chromatic number at most $\tau_1$.
Let $\mathcal{C}_2$ be the ideal of $G\in \mathcal{C}$ with clique number at most $\kappa$ and no $d$-peripheral $\ell$-hole. We suppose
that there are graphs in $\mathcal{C}_2$ with arbitrarily large chromatic number, and so $\mathcal{C}_2$ is not 2-controlled,
by \ref{rad2}. Consequently there exists $\tau_2$ such that if $\mathcal{C}_3$ denotes the class of graphs $G\in \mathcal{C}_2$
with $\chi^2(G)\le \tau_2$, there are graphs in $\mathcal{C}_3$ with arbitrarily large chromatic number. Hence by \ref{manyholes}
with $\rho=3$, $\mathcal{C}_3$ is not 3-controlled, and so on; and we deduce that there is an ideal $\mathcal{C}_{\rho}$ of graphs
in $\mathcal{C}$ that is not $\rho$-controlled, a contradiction. This proves \ref{radthm2}.~\bbox

\section{Controlling 8-balls}\label{sec:uncontrolled}

In this section we prove \ref{uncontrolled}.
We use the following relative of \ref{greentouchrad}, proved in~\cite{longoddholes}:

\begin{thm}\label{greentouch}
Let $\tau,c\ge 0$ and let $(W_1\ll W_n)$ be a $\tau$-colourable grading of a graph $G$. Let $H$ be a subgraph of $G$ (not necessarily
induced)
with $\chi(H)>\tau+2(c+\chi^1(G))$. Then there is an edge $uv$ of $H$, and a subset $X$ of $V(G)$,
such that
\begin{itemize}
\item $G[X]$ is connected;
\item $u,v$ are both earlier than every vertex in $X$;
\item $v$ has a neighbour in $X$, and $u$ does not; and
\item $\chi(X)>c$.
\end{itemize}
\end{thm}
We deduce a version of \ref{greenedgerad} that has no assumption of $\rho$-control:

\begin{thm}\label{greenedge}
Let $G$ be a graph, and let $B,C\subseteq V(G)$, where $B$ covers $C$. Let every induced subgraph $J$ of $G$ with
$\omega(J)<\omega(G)$
have chromatic number at most $\tau$.
Let the enumeration $(b_1\ll b_m)$ of $B$ and the grading $(W_1\ll W_n)$ of $G[C]$ be compatible. 
Let $(W_1\ll W_n)$ be $\tau$-colourable, and let $\chi(G[C])>\tau^2 (2c+3\tau)$.
Then there is a square edge $uv$, and a subset $X$ of $V(G)$,
such that
\begin{itemize}
\item $G[X]$ is connected;
\item $u,v$ are both earlier than every vertex in $X$;
\item $v$ has a neighbour in $X$, and $u$ does not; and
\item $\chi(X)>c$.
\end{itemize}
\end{thm}
\Proof Let $H$ be as in \ref{getgreenedge}; then by \ref{getgreenedge}, $\chi(G[C])\le \tau^2\chi(H)$.
Since $\chi(G[C])>\tau^2 (2c+3\tau)$ and $\chi^1(G)\le \tau$, it follows that $\chi(H)>\tau+2(c+\chi^1(G))$.
By \ref{greentouch} applied to $G[C]$ and $H$, we deduce
that there is an edge $uv$ of $H$, and a subset $X$ of $V(G)$, satisfying the theorem. This proves \ref{greenedge}.~\bbox

A {\em shower} in $G$ is a sequence $(L_0, L_1\ll L_k,s)$ where $L_0, L_1\ll L_k$ are pairwise disjoint subsets
of $V(G)$ and $s\in L_k$,
such that
\begin{itemize}
\item $|L_0|=1$;
\item $L_{i-1}$ covers $L_i$ for $1\le i< k$;
\item for $0\le i<j\le k$, if $j>i+1$ then no vertex in $L_j$ has a neighbour in $L_i$; and
\item $G[L_k]$ is connected.
\end{itemize}
(Note that we do not require that $L_{k-1}$ covers $L_k$.)
We call the vertex in $L_0$ the {\em head} of the shower, and
$s$ its {\em drain}, and $L_0\cup \cdots\cup L_k$ is its {\em vertex set}. 
(Thus the drain can be any vertex of $L_k$.)
The set of vertices in $L_k$ with a neighbour in $L_{k-1}$
is called the {\em floor} of the shower. 

Let $\mathcal{S}$ be a shower with head $z_0$, drain $s$ and vertex set $V$.
An induced path of $G[V]$ between $z_0,s$ is called a {\em jet}
of $\mathcal{S}$. For $d\ge 0$, a jet $J$ is {\em $d$-peripheral} if there is a subset $X$ of the floor of the shower, anticomplete to $V(J)$,
with $\chi(G[X])>d$.
The set of all lengths of $d$-peripheral jets of $\mathcal{S}$ is called the {\em $d$-jetset} of $\mathcal{S}$.
For integers $\ell\ge 2$ and $1\le k\le \ell$, we say a $d$-jetset is {\em $(k,\ell)$-complete} if there are
$k$ jets $J_0\ll J_{k-1}$, all $d$-peripheral, such that $|E(J_j)|=|E(J_0)|+j$ modulo $\ell$ for $0\le j\le k-1$.

We will prove:

\begin{thm}\label{multijet}
Let $\tau\ge 0$ and $\ell\ge 2$. For all integers $d\ge 0$ and $t$ with $1\le t\le \ell$ there exists $c_{d,t}\ge 0$ with the following property.
Let $G$ be a graph such that $\chi^{3}(G)\le \tau$, and such that
every induced subgraph $J$ of $G$ with
$\omega(J)<\omega(G)$
has chromatic number at most $\tau$. 
Let $\mathcal{S}=(L_0, L_1\ll L_k,s)$ be a shower in $G$,
with floor of chromatic number more than $c_t$. Then the $d$-jetset of $\mathcal{S}$ is $(t,\ell)$-complete.
\end{thm}
\Proof Suppose first that $t=1$. Let $c_{d,1}=d+\tau$, and let $G$ and 
$\mathcal{S}=(L_0, L_1\ll L_k,s)$ be as in the theorem, with floor $F$ where
$\chi(F)>c_1$. Since $c_{d,1}\ge 0$, it follows that $F\ne \emptyset$, and so there is an induced path $P$ of $G[L_k]$ 
between $s$ and some vertex in $F$. Choose such a path $P$ of minimum length, and let its end in $F$ be $v$. If $s=v$ let 
$u=v$, and otherwise let $u$ be the neighbour of $v$ in $P$. It follows that no vertex of $P$ belongs to or has a neighbour in $F$
except $u,v$.
Also there is a path $Q$ of length $k$ between $v$ and $z_0\in L_0$, with 
one vertex $w\in L_{k-1}$ and all others in $L_0\cup\cdots\cup L_{k-2}\cup \{v\}$. Thus $P\cup Q$ is a jet. Moreover, the only 
vertices of $P\cup Q$ that belong to or have a neighbour in $F$ are $u,v,w$, and so all these neighbours have
$G$-distance at most two from $v$, and consequently have chromatic number at most $\tau$.
Since $\chi(F)>c_{d,1}=d+\tau$, it follows that $P\cup Q$ is a $d$-peripheral jet, as required.

We may therefore assume that $2\le t\le \ell$, and inductively the result holds for $t-1$ (and all $d$).
Define $d'=d+2\tau$, and let $c_{d,t}=2(\ell+1)\tau^2 (2c_{d',t-1}+7\tau)$; let $G$ be a graph, and let $\mathcal{S}=(L_0, L_1\ll L_k,s)$ be a shower in $G$
with floor $F$ of chromatic number more than $c_{d,t}$.
For $i\ge 0$, let $M_i$ be the set of vertices with $G[L_k]$-distance exactly $i$ from $s$. Then there exists $r\ge 0$ such that
$\chi(F\cap M_r)\ge \chi(F)/2$; and $r\ge 3$ since $\chi^2(G)\le \tau<\chi(F)/2$.
For $0\le i\le r$, each vertex $v\in M_i$ is joined to $s$ by an induced path of length $i$ with interior in
$M_1\cup\cdots\cup M_{i-1}$; let us call such a path a {\em bloodline} of $v$.
\begin{figure}[H]
\centering
\tikzset{snake it/.style={decorate, decoration=snake}}
\begin{tikzpicture}[scale=.5,auto=left]
\draw (0,0) ellipse (10cm and 1cm);
\draw (0,3) ellipse (10cm and 1cm);
\node at (5,0) {$L_{k-1}$};
\node at (5,3) {$L_{k-2}$};
\node at (5,-3) {$L_k$};
\node at (-9,-6.5) {$M_r$};
\node at (-5,-6.5) {$M_{r-1}$};
\node at (-3,-6.5) {$M_{r-2}$};
\node at (9.4,-6.5) {$M_0$};
\node [left] at (-8,-4.5) {$u$};
\node [right] at (-1,-3.5) {$v$};
\node [left] at (-5,0) {$b$};

\draw [rounded corners] (-10,-2) rectangle (10,-7);
\draw [dashed] (0,4.2) -- (0,5.5);
\draw [dashed] (-1.0,-4.5) -- (6.7,-4.5);
\draw [very thin] (0,2.3) -- (0,0.7);
\draw [very thin] (-7,2.5) -- (-7,0.5);
\draw [very thin] (7,2.5) -- (7,0.5);
\draw [very thin] (0,-0.6) -- (0,-2.5);
\draw [very thin] (-7,-0.5) -- (-7,-2.5);
\draw [very thin] (7,-0.5) -- (7,-2.5);
\node [above] at (9.5,-4.5) {$s$};
\draw [very thin] (8.8,-2) -- (8.8,-7);
\draw [very thin] (7.6,-2) -- (7.6,-7);
\draw [very thin] (-2,-2) -- (-2,-7);
\draw [very thin] (-4,-2) -- (-4,-7);
\draw [very thin] (-6,-2) -- (-6,-7);
\tikzstyle{every node}=[inner sep=1.5pt, fill=black,circle,draw]
\node (s) at (9.5,-4.5) {};
\node (u) at (-8,-4.5) {};
\node (v) at (-1,-3.5) {};
\node (b) at (-5,0) {};
\draw [very thick] (u) -- (b);
\draw [very thick] (v) -- (b);

\end{tikzpicture}

\caption{$v$ is a grandparent of $u$.} \label{fig:3}
\end{figure}
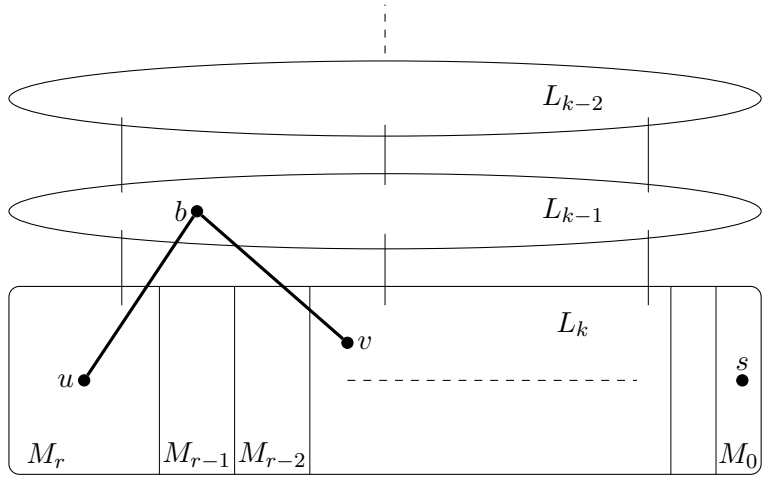

If $u\in F\cap M_r$ and $v\in M_0\cup\cdots\cup M_{r-2}$, we say that $v$ is a 
{\em grandparent} of $u$ if there exists $b\in L_{k-1}$ adjacent to both $u,v$. 
Let $C$ be the set of vertices in $M_r$ with a grandparent in $M_0\cup\cdots\cup M_{r-2}$.
\\
\\
(1) {\em If $\chi(C)>\ell \tau^2 (2c_{d',t-1}+7\tau)$ then the theorem holds.}
\\
\\
Let $(v_1\ll v_n)$ be an enumeration of $M_0\cup\cdots \cup M_{r-2}$, where the vertex in $M_0$ is first, followed by the vertices
in $M_1$ in some order, and so on; more precisely, for $1\le i<j\le n$, if $v_i\in M_a$
and $v_j\in M_b$ where $a,b\in \{0\ll r-2\}$, then $a\le b$. Let $B_0$ be the set of vertices
in $L_{k-1}$ that have neighbours in $M_0\cup\cdots \cup M_{r-2}$, and let $(b_1\ll b_m)$ be an enumeration of $B_0$,
enumerating the members
of $B_0$ with the earliest neighbours in $(v_1\ll v_n)$ first; more precisely,
for $1\le i<j\le m$, if $b_j$ is adjacent to $v_q$ for some $q\in \{1\ll n\}$, then there exists $p\in \{1\ll q\}$
such that $b_i$ is adjacent to $v_p$. 
We say $v_i$ is the {\em earliest grandparent}
of $u\in M_r$ if $i$ is minimum such that $v_i$ is a grandparent of $u$.
For $1\le i\le n$, let $W_i$ be the set of vertices in $M_r$ whose earliest grandparent is $v_i$.
Thus $C=W_1\cup\cdots\cup W_n$, and $(W_1\ll W_n)$ is a grading of
$G[C]$. It is $\chi^2(G)$-colourable, since $W_i\subseteq N^2_G[v_i]$ for each $i$, and hence $\tau$-colourable. 
For $u\in C$, if $v$ is the earliest grandparent of $u$, then the earliest parent of $u$ is adjacent to $v$. Consequently
the enumeration $(b_1\ll b_m)$ and the
grading $(W_1\ll W_n)$ are compatible; because if $u,v\in W$ with $u$ earlier than $v$, the earliest grandparent of $u$ is
earlier than
the earliest grandparent of $v$, and consequently the earliest parent of $u$ is
earlier than
the earliest parent of $v$.

For $0\le j<\ell$, let $C_j$ be the set of vertices $u\in C$ such that $i-j$ is a multiple of $\ell$, where $i$ is the length
of the bloodline of the earliest grandparent of $u$. Thus $C_0\cup\cdots\cup C_{\ell-1} = C$.
Choose $j$ such that $\chi(C_j)>\tau^2 (2c_{d,t-1}+7\tau)$.  Then by \ref{greenedge}, 
there is a square edge $uv$ with $u,v\in C_j$, and a subset $X$ of $C_j$,
such that
\begin{itemize}
\item $G[X]$ is connected;
\item $u,v$ are both earlier than every vertex in $X$;
\item $v$ has a $G$-neighbour in $X$, and $u$ does not; and
\item $\chi(X)>c_{d,t-1}+2\tau$.
\end{itemize}

\begin{figure}[H]
\centering
\begin{tikzpicture}[scale=.5,auto=left]
\draw (10,1) arc (90:270:20cm and 1cm);
\node at (5,0) {$L_{k-1}$};
\node at (5,-7) {$L_k$};
\node at (-8,-11) {$M_r$};
\node at (-5,-11) {$M_{r-1}$};
\node at (-3,-11) {$M_{r-2}$};
\node at (-9,-4) {$X$};
\node at (-9,-8.5) {$W_i$};
\node at (-9,-7.5) {$W_{i+1}$};
\node [left] at (-7,-8.5) {$v$};
\node [right] at (4,-4.5) {$v''$};
\node [left] at (1,0) {$v'$};
\node [left] at (-6.5,-6.5) {$u$};
\node [right] at (1,-5.5) {$u''$};
\node [left] at (-1,0) {$u'$};

\draw [rounded corners] (10,-12) -- (-10,-12) -- (-10,-2) -- (10,-2);
\draw [dashed] (0,1.2) -- (0,2.5);
\draw [dashed] (-1.5,-7) -- (0,-7);
\draw [very thin] (0,-0.6) -- (0,-2.5);
\draw [very thin] (-7,-0.2) -- (-7,-2.3);
\draw [very thin] (7,-0.6) -- (7,-2.5);
\draw [very thin] (-2,-2) -- (-2,-12);
\draw [very thin] (-4,-2) -- (-4,-12);
\draw [very thin] (-6,-2) -- (-6,-12);
\draw [very thin] (-10,-9) -- (-6,-9);
\draw [very thin] (-10,-8) -- (-6,-8);
\draw [very thin] (-10,-7) -- (-6,-7);
\draw [very thin] (-10,-6) -- (-6,-6);
\draw [very thin] (-10,-5) -- (-6,-5);
\draw [dashed] (-8,-9.2) -- (-8,-10);
\draw [dashed] (-8,-4.8) -- (-8,-4);
\draw (-9,-4) ellipse (0.7cm and 1.4cm);

\tikzstyle{every node}=[inner sep=1.5pt, fill=black,circle,draw]
\node (v) at (-7,-8.5) {};
\node (v'') at (4,-4.5) {};
\node (v') at (1,0) {};
\node (u) at (-6.5,-6.5) {};
\node (x) at (-9,-4.7) {};
\node (u'') at (1,-5.5) {};
\node (u') at (-1,0) {};

\draw [very thick] (v) -- (v');
\draw [very thick] (v') -- (v'');
\draw [very thick] (u) -- (v);
\draw [very thick] (x) -- (v);
\draw [very thick] (u) -- (u');
\draw [very thick] (u') -- (u'');

\end{tikzpicture}

\caption{The square edge and its relatives.} \label{fig:4}
\end{figure}
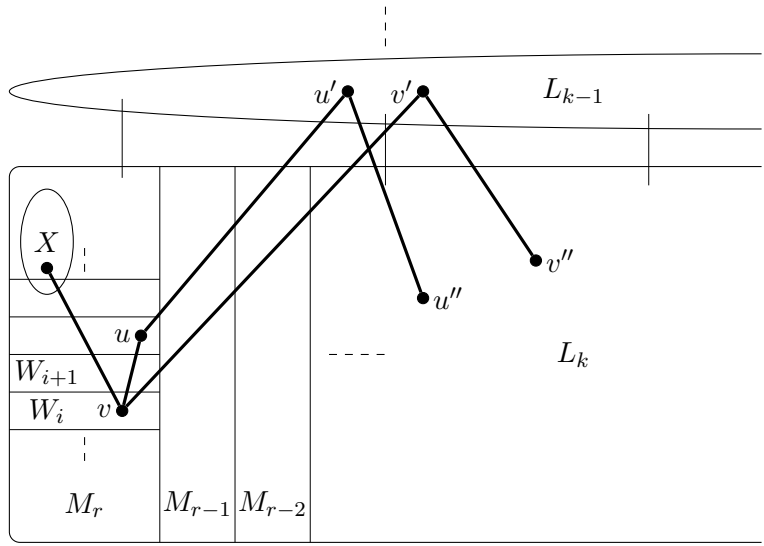

Let the earliest parents of $u,v$ be $u', v'$ respectively, and let their earliest grandparents be $u'', v''$ respectively.
Let $P$ be the induced path between $v,s$ consisting of the path $v\d v'\d v''$ and
a bloodline of $v''$.
Thus $P$ has length $j+2$ modulo $\ell$. Let $Q$ be the path between $v,s$
consisting of the edge $uv$, the path $u\d u'\d u''$, and a bloodline of $u''$. Note that $Q$ is induced,
since $v$ is not adjacent to $u'$ (because $uv$ is square). Moreover, $Q$ has length $j+3$ modulo $\ell$.

Let $Z$ be the set of vertices in $X$ with $G$-distance at least four from both $u',v'$, and let $L_{k-1}'$ be the set of 
vertices in $L_{k-1}$ with a neighbour in $Z$. Let $L_{k-2}'$ be the set of vertices in $L_{k-2}$ with a 
neighbour in $L_{k-1}'$, and for $0\le i\le k-3$
let $L_i'=L_i$. It follows that $L_{i-1}$ covers $L_i$ for $0\le i\le k-2$, and so $(L_0'\ll L_{k-1}',X\cup \{v\},v)$ is a 
shower $\mathcal{S}'$. Let $V'$ be its vertex set.

We claim that $V'\cap V(P\cup Q)=\{v\}$, and every edge between $V'$ and $V(P\cup Q)$ is incident with $v$. For suppose that 
$a\in V'$ and $b\in V(P\cap Q)$ are equal or adjacent, and $a,b\ne v$. If $b=u$, then $a\ne b$ since $u\notin V'$; $a\notin X$
since $u$ has no neighbour in $X$; and so $a\in L_{i}'$ for some $i'<k$. Then $i=k-1$ since $b\in L_k$, and so $a$ has a neighbour 
in $Z$ from the definition of $L_{k-1}'$; and so the $G$-distance between $a,u$ is at least two from the definition of $Z$,
a contradiction. Thus $b\ne u$. Next suppose that $b\in \{u',v'\}$.
Since $u',v'$ are the earliest parents of $u,v$ respectively,
and $u,v$ are earlier than every vertex in $X$, it follows that no vertex in $X$ is adjacent to $b$; and so $a\ne b$, and $a\notin X$.
Hence $a\in L_i'$ for some $i<k$, and since $b\in L_{k-1}$ it follows that $k-2\le i\le k-1$. But then the $G$-distance between $a$
and some vertex of $Z$ is at most two, and since the $G$-distance between $b,Z$ is at least four, it follows that $a,b$ are nonadjacent,
a contradiction. So $b\notin \{u',v'\}$ and so $b\in M_j$ for some $j\le r-2$, and so $b$ is a vertex of a bloodline of the earliest
grandparent of one of $u,v$ (say $u'', v''$ respectively). Consequently $a\ne b$, and $a\notin X$, and so 
$a\in L_{k-1}'$. Choose $z\in Z$ adjacent to $a$; then $u,v$ are earlier than $z$, and so $u'',v''$ are both
earlier than the earliest grandparent of $z$. It follows that no parent of $z$ has a neighbour in a bloodline of either of $u'',v''$,
and so $a,b$ are nonadjacent, 
a contradiction. This proves our claim that $V'\cap V(P\cup Q)=\{v\}$, and every edge between $V'$ and $V(P\cup Q)$ is incident with $v$.
In particular $Z$ is anticomplete to $V(P\cup Q)$.

The floor of $\mathcal{S}'$
includes $Z$, and $\chi(Z)\ge \chi(X)-2\chi^3(G)>c_{d',t-1}$. Consequently the $d'$-jetset of $\mathcal{S}'$ is $(t-1,\ell)$-complete;
let $J_0\ll J_{t-2}$ be corresponding $d'$-peripheral jets of $\mathcal{S}'$. For $0\le h\le t-2$, both $J_h\cup P$ and $J_h\cup Q$
are jets of $\mathcal{S}$, and we claim they are $d$-peripheral. For let $0\le h\le t-2$, and let $Y$ be a subset of the floor of
$\mathcal{S}'$ with $\chi(G[Y])>d'$, anticomplete to $V(J_h)$. Thus $Y\subseteq X\cup \{v\}\subseteq F$. Let $Y'=Y\cup Z$.
Since every vertex in $Y\setminus Z$ has $G$-distance at most three from one of $u',v'$, it follows that
$\chi(Y\setminus Z)\le 2\tau$, and so $\chi(Y')\ge \chi(Y)-2\tau> d$. Since $Y'$ is anticomplete to $V(P\cup Q)$, 
this proves our claim that both $J_h\cup P$ and $J_h\cup Q$ are $d$-peripheral jets of $\mathcal{S}$. Consequently the $d$-jetset of 
$\mathcal{S}$ is $(t,\ell)$-complete. This proves (1).
\\
\\
(2) {\em If $\chi((F\cap M_r)\setminus C)>\tau^2 (2c_{d',t-1}+7\tau)$ then the theorem holds.}
\\
\\
Let us write $C' = (F\cap M_r)\setminus C$; then $C'$ is the set of vertices in $M_r$ that have neighbours in $L_{k-1}$,
but every such neighbour has no neighbour in $M_0\cup\cdots\cup M_{r-2}$. The neighbours in $L_{k-1}$ might or
might not have neighbours in $M_{r-1}$. Take an arbitrary enumeration $(b_1\ll b_n)$ of $M_{r-1}$, and for $1\le i\le n$
let $W_i$ be the set of vertices in $C'$ adjacent to $b_i$ and nonadjacent to $b_1\ll b_{i-1}$. Thus $(W_1\ll W_n)$ is a grading
of $G[C']$, compatible with $(b_1\ll b_n)$, and it is $\chi^1(G)$-colourable and hence $\tau$-colourable. By \ref{greenedge},
there is a square edge $uv$ of $G[C']$, and a subset $X$ of $C'$,
such that
\begin{itemize}
\item $G[X]$ is connected;
\item $u,v$ are both earlier than every vertex in $X$;
\item $v$ has a neighbour in $X$, and $u$ does not; and
\item $\chi(X)>c_{d',t-1}+2\tau$.
\end{itemize}
Let $u',v'$ be the earliest neighbours in $M_{r-1}$ of $u,v$ respectively.
Let $P$ be an induced path between $v,s$ consisting of the edge $vv'$ and a bloodline of $v'$; and let $Q$ be the path
consisting of the edges $uv,uu'$, and a bloodline of $u'$. They are both induced. 
Let $Z$ be the set of all vertices in $X$ with $G$-distance at least four from both of $u,v$. Let $L_{k-1}'$
be the set of vertices in $L_{k-1}$ with a neighbour in $Z$, and for $0\le i\le k-2$ let $L_i'=L_i$. Then 
$(L_0'\ll L_{k-1}', X\cup \{v\},v)$ is a shower $\mathcal{S}'$. Moreover, its vertex set $V'$ satisfies $V'\cap V(P\cup Q)=\{v\}$,
and every edge between $V'$ and $V(P\cup Q)$ is incident with $v$ (the proof is as in (1), and we omit it).
Its floor includes $Z$, and $\chi(Z)\ge \chi(X)-2\chi^3(G)>c_{d',t-1}$;
and so the $d'$-jetset of $\mathcal{S}'$ is $(t-1,\ell)$-complete. But for each jet $J$ of $\mathcal{S}'$, $J\cup P, J\cup Q$
are both jets of $\mathcal{S}$, and as in (1) it follows that the $d$-jetset of $\mathcal{S}$ is $(t,\ell)$-complete. This proves (2).

\bigskip

Since $\chi(F)>c_{d,t}$, it follows that $\chi(F\cap M_r)>(\ell+1)\tau^2 (2c_{d',t-1}+7\tau)$, and so
either $\chi(C)>\ell \tau^2 (2c_{d',t-1}+7\tau)$ or $\chi((F\cap M_r)\setminus C)>\tau^2 (2c_{d',t-1}+7\tau)$. Hence
the result follows from (1) or (2). This proves \ref{multijet}.~\bbox

A {\em recirculator} for a shower $(L_0, L_1\ll L_k,s)$ with head $z_0$ is an induced path $R$
with ends $s,z_0$ such that no internal vertex of $R$ belongs to $V$
and no internal vertex of $R$ has any neighbours in $V\setminus \{s,z_0\}$.
We need the following, proved in~\cite{holeseq}:

\begin{thm}\label{doubleshower}
Let $c\ge 0$ be an integer, and let $G$ be a graph such that $\chi(G)>44c+4\chi^{8}(G)$.
Then there is a shower in $G$, with floor of chromatic number more than $c$, and with a recirculator.
\end{thm}

We deduce \ref{uncontrolled}, which we restate:

\begin{thm}\label{uncontrolled2}
For all integers $\ell\ge 2$ and $\tau,d\ge 0$ there is an integer $c\ge 0$ with the following property.
Let $G$ be a graph such that $\chi^8(G)\le \tau$, and every induced subgraph $J$ of $G$ with
$\omega(J)<\omega(G)$
has chromatic number at most $\tau$.
If $\chi(G)>c$ then there are $\ell$ $d$-peripheral holes in $G$ with lengths of all possible values modulo $\ell$.
\end{thm}
\Proof
Let $c_{d,\ell}$ be as in \ref{multijet}, and let $c=44c_{d,\ell}+4\tau$. Let $G$ be a graph 
such that $\chi^8(G)\le \tau$, and every induced subgraph $J$ of $G$ with
$\omega(J)<\omega(G)$
has chromatic number at most $\tau$, with $\chi(G)>c$. By \ref{doubleshower} there
is a shower in $G$, with floor of chromatic number more than $c_{d,\ell}$, and with a recirculator.
By \ref{multijet} the $d$-jetset of this shower is $(\ell,\ell)$-complete. Thus adding the recirculator to each of the corresponding jets
gives the $\ell$ $d$-peripheral holes we need. This proves \ref{uncontrolled2}.~\bbox

\section{Two consecutive holes}\label{sec:2holes}

Finally let us prove \ref{2holes}, which we restate.

\begin{thm}\label{2holesagain}
For each $\kappa,\ell\ge 0$ there exists $c\ge 0$ such that every graph $G$ with $\omega(G)\le \kappa$ and $\chi(C)>c$
has holes of two consecutive lengths, both of length more than $\ell$.
\end{thm}
\Proof
We may assume that $\ell\ge 8$. 
By induction on $\kappa$, there exists $\tau_1$ such that 
every graph with clique number less than $\kappa$ and chromatic number
more than $\tau_1$ has two holes of consecutive length more than $\ell$. 
Let $\mathcal{C}_2$ be the ideal of graphs with clique number at most $\kappa$
and with no two holes of consecutive lengths more than $\ell$. By \ref{rad2}, $\mathcal{C}_2$ is not 2-controlled, and so 
for some $\tau_1$ there are graphs $G$ in $\mathcal{C}_2$ with arbitrarily large chromatic number and $\chi^2(G)\le \tau$.
Let $\mathcal{C}_3$ be the ideal of graphs in $G$ with $\chi^2(G)\le \tau$. By \ref{manyholes} with $\rho=3$,
$\mathcal{C}_3$ is not 3-controlled, and so on; and we deduce that there is an ideal $\mathcal{C}_{\ell}$ of graphs in $\mathcal{C}_2$,
with unbounded chromatic number, and a number $\tau$ such that $\chi^{\ell}(G)\le \tau$ for each $G\in \mathcal{C}_{\ell}$.

Let $c'=14\tau^3$, and let $c=44c'+4\tau$; and choose $G\in \mathcal{C}_{\ell}$ with $\chi(G)>c$. 
By \ref{doubleshower}, there is a shower $\mathcal{S}=(L_0\ll L_k,s)$ in $G$ with a recirculator and with floor $F$ with chromatic number
more than $14\tau^3$. Define $d, M_0\ll M_d$ and $C$
as in the proof of \ref{multijet}.
\\
\\
(1) {\em If $\chi(C)> 7\tau^3$ then the theorem holds.}
\\
\\
Define the enumeration $(v_1\ll v_n)$ of $M_0\cup \cdots\cup M_{d-2}$,
and $B_0$ and its enumeration $(b_1\ll b_m)$,
as in the proof of \ref{multijet}. Let $b_{m+1-i}'=b_i$; so $(b_1'\ll b_m')=(b_m\ll b_1)$ is also an enumeration of the same set.
For $i=1\ll m$ let $W_i$ be the set of vertices in $M_d$ that are adjacent to $b_i'$ and to none of $b_1'\ll b_{i-1}'$. Thus
$W_1\cup\cdots\cup W_m=C$, and $(W_1\ll W_m)$ is a $\chi^1(G)$-colourable (and hence $\tau$-colourable) grading of $G[C]$,
compatible with $(b_1'\ll b_m')$. By \ref{greenedge} (taking $c=2\tau$),
there is a square edge $uv$ of $G[C]$, and a subset $X$ of $C$,
such that
\begin{itemize}
\item $G[X]$ is connected;
\item $u,v$ are both earlier than every vertex in $X$;
\item $v$ has a neighbour in $X$, and $u$ does not; and
\item $\chi(X)>2\tau$.
\end{itemize}
Let $u',v'\in B_0$ be the earliest parents of $u,v$ respectively. Since $\chi(X)>2\tau$, there exists $x\in X$ with $G$-distance
at least four from each of $u',v'$. Let $x'\in B_0$ be its earliest parent, and let $R$ be an induced path of $G[X\cup \{v,x'\}]$ 
between $v$ and $x'$.
Now no vertex of the interior of $R$ is adjacent to either of $u',v'$ since $u,v$ are earlier than every member of $X$
and $u',v'$ are their earliest parents. Also, $x'$ is nonadjacent to $u,v,u',v'$ since the $G$-distance between $x$ and $u',v'$
is at least four. Since $uv$ is square, the path $Q'$ obtained from $R$ by adding the edge $vv'$ is induced,
and since $u$ has no neighbour in $X$, so is the path $P'$ obtained from $R$ by adding the path $v\d u\d u'$. 
Now there are paths between the apex of $\mathcal{S}$ (say $a$)
and $u'$, and between $a,v'$, both of length $k-1\ge \ell$. No vertex of $L_k$ has a neighbour different from $u',v'$ in these paths.
Also $x'$ has no neighbour in $P'\cup Q'$; because any such neighbour would be in $L_k\cup L_{k-1}\cup L_{k-2}$, and hence
would have $G$-distance at most one from one of $u',v'$, which is impossible since the $G$-distance between $x$ and $u',v'$
is at least four. Thus by taking the union of the first of these with $P'$ and the second with $Q'$,
we obtain two paths $P,Q$, both induced and both between $a,x'$, with consecutive lengths. Choose $i$ minimum such that $x'$
is adjacent to $v_i$, and let $R$ be a bloodline of $v_i$ (defined as in \ref{multijet}). Let $u'=b_f'$, $v'=b_g'$, $x'=b_h'$. Since $u,v$ are earlier than $x$,
it follows that $f,g<h$. Now $u',v'$ are nonadjacent to $v_i$ since the $G$-distance between $u'v'$ and $x$ is at least four.
Since $u'=b_{m+1-f}$ and $x'=b_{m+1-h}$, and $m+1-f>m+1-h$, it follows that no neighbour of $u'$ belongs to $V(R)$, and similarly
for $v'$. Thus the union of $P$ with the edge $x'v_i$ and $R$ is induced, and so is the union of $Q$ with $x'v_i$ and $R$.
Consequently there are jets of $\mathcal{S}$ with consecutive lengths. By taking their unions with the recirculator, we obtain
holes of consecutive lengths. This proves (1).
\\
\\
(2) {\em If $\chi((F\cap M_d)\setminus C)>7\tau^3$ then the theorem holds.}
\\
\\
The proof of this is the same as that for step (2) of the proof of \ref{multijet} and we omit it.

\bigskip
From (1) and (2) the result follows.~\bbox

\section{Some connections with homology}

In the 1990s, Kalai and Meshulam made several intriguing conjectures connecting the chromatic number of a graph with the homology 
of a simplicial complex associated with $G$. (Most of them are mentioned in~\cite{kalai}, and see also~\cite{kalai2}.)  
The {\em $n$th Betti number}
$b_n(X)$ of 
a simplicial complex $X$ is the rank of the $n$th homology group $H_n(X)$ (see, for instance, \cite{hatcher}).  The {\em Euler characteristic\footnote{The 
standard notation for the Euler characteristic of a simplicial complex $X$ is $\chi(X)$; however, we will avoid using that 
notation here, as there is clearly some potential for confusion between $\chi(I(G))$ and $\chi(G)$.}} of $X$ is 
$\sum_{n\ge 0}(-1)^n c_n(X)$, where $c_n(X)$ is the number of $n$-faces in $X$; it turns out that the Euler characteristic is 
also equal to the alternating sum $\sum_{n\ge0}(-1)^n b_n(X)$ of Betti numbers.
The {\em independence complex} $I(G)$ of a graph $G$ is the simplicial complex whose faces are the stable sets of vertices.  For a graph $G$, we say that its {\em total Betti number} is $b(G):=\sum_{n\ge 0}b_n(G)$.  The total Betti number of a graph $G$ is clearly greater than or equal to the modulus of the Euler characteristic of $I(G)$, as the former is the sum of Betti numbers and the latter is equal to the alternating sum.

Kalai and Meshulam made several conjectures on this topic: one was
already mentioned (at \ref{bonamy}), and two others are as follows:

\begin{thm}\label{kalaiconj2}
For every integer $k\ge 0$ there exists $c$  such that the following holds.
If $b(H)\le k$ for every induced subgraph $H$ of $G$ then $\chi(G)\le c$.
\end{thm}

\begin{thm}\label{kalaiconj1}
For every integer $k\ge 0$ there exists $c$  such that following holds.  For every graph $G$, if the 
Euler characteristic of $I(H)$ has modulus at most $k$ for every induced subgraph $H$ of $G$ then 
$\chi(G)\le c$.
\end{thm}

Kalai and Meshulam also asked about the graphs $G$ that satisfy the condition in \ref{kalaiconj1} with $k=1$ (i.e.~for every induced subgraph $H$, the Euler characteristic of $H$ lies in $\{-1,0,1\}$).  They conjectured that $G$ has this property if and only if 
$G$ has no induced cycle of length divisible by three.  We prove this conjecture in~\cite{ternary}, with Chudnovsky and Spirkl.  

In this paper, we prove conjectures \ref{kalaiconj2} and \ref{kalaiconj1}.
The second conjecture is clearly stronger, as the modulus of the Euler characteristic of $I(G)$ is at most the total Betti number of $G$.
In this section, we will prove both of these conjectures, using \ref{superkalai}.

Say that
a graph $G$ is {\em $k$-balanced} if for every induced subgraph $H$ of $G$, the number of stable sets in $H$ of even cardinality differs by at most $k$ from the
number of stable sets of odd cardinality.  The condition in \ref{kalaiconj1} is exactly that $G$ is $k$-balanced, so \ref{kalaiconj1} (and therefore also \ref{kalaiconj2})  is an immediate consequence of the following result.
\begin{thm}\label{kalaiconj}
For every integer $k\ge 0$ there exists $c$  such that $\chi(G)\le c$ for every $k$-balanced graph $G$.
\end{thm}
\noindent{\bf Proof of \ref{kalaiconj}, assuming \ref{superkalai}.\ \ }
Let $k\ge 1$ be an integer. By \ref{superkalai}, we may choose $c$ such that every graph $G$ with $\omega(G)\le k+1$
and $\chi(G)>c$ contains $k$ holes, pairwise anticomplete  and each of length a multiple of three. We claim that
every $k$-balanced graph has chromatic number at most $c$. For if $G$ is $k$-balanced, then $G$ has no clique of cardinality more
than $k+1$ (because a complete subgraph $H$ with $k+2$ vertices has $k+2$ odd stable sets and only one even one); and $G$ does not
have $k$ holes that are pairwise anticomplete, each of length a multiple of three. (We leave the reader to check this.) This
proves \ref{kalaiconj}.~\bbox


\section*{Acknowledgement}
We would like to thank Gil Kalai for helpful discussions.

\end{document}